\documentstyle[12pt,amscd]{amsart}
\def\cite#1{[#1]}
\newtheorem{theorem}{Theorem}[section]
\newtheorem{proposition}[theorem]{Proposition}
\newtheorem{corollary}[theorem]{Corollary}
\newtheorem{lemma}[theorem]{Lemma}
\theoremstyle{definition}
\newtheorem{definition}[theorem]{Definition}
\newtheorem{example}[theorem]{Example}
\newtheorem{remark}[theorem]{Remark}
\newtheorem{problem}[theorem]{Problem}
\newenvironment{proof}{\trivlist \item[\hskip \labelsep{\sc
Proof.\kern1pt}]}{\endtrivlist}
\newenvironment{sketch}{\trivlist \item[\hskip \labelsep{\sc
Sketch.\kern1pt}]}{\endtrivlist}
\newenvironment{proofnodot}{\trivlist \item[\hskip \labelsep{\sc
Proof}]}{\endtrivlist}

\def\th#1{\xmode{{#1}^{\operatoratfont th}}}
\def\squareSE#1#2#3#4{\diagramx{#1&\mapE{}&#2\cr
		      \mapS{}&&\mapS{}\cr #3&\mapE{}&#4\cr}}
\def\geom{g\'eom\'etrie}
\def\setofh#1{\{\hbox{#1}\}}
\def\dmapx[[#1||#2]]{%
     \setbox0=\hbox{$\displaystyle #1\ \ \mapE\ \ \ #2$}%
     \ifdim\wd0<\hsize$$#1\ \ \mapE\ \ \ #2$$\else%
     \begin{flushleft} $\displaystyle #1$ \end{flushleft}
     \begin{flushright} $\displaystyle\mapE{}\ \ \ #2$ \end{flushright}\fi}
\def\NN{\Bbb N}
\def\CP{{\Bbb C}\kern1pt{\Bbb P}}
\def\tC{\tilde C}
\def\CM{Cohen-Macaulay}
\def\vec#1#2#3{#1_{#2}, \ldots, #1_{#3}}
\def\SPEC{\mathop{\mathbf{Spec}}\nolimits}

\def\oG{\overline{G}}
\def\etale{\'etale}
\def\pref#1{(\ref{#1})}
\def\iso{\cong}
\def\IN{\subset}
\def\et{_{\hbox{\footnotesize\'et}}}
\def\FG{finitely generated}
\def\WMAT{we may assume that}

\def\includeE#1{{\lhook\kern-3.5pt\joinrel\smash{
    \mathop{\longrightarrow}\limits^{#1}}}}
\def\includeS#1{\hbox{\raise11.0pt\hbox{$\scriptscriptstyle\cap$}\kern-4.3pt
    \lower2.5pt\hbox{${\Big\downarrow}$}}
    \rlap{$\vcenter{\hbox{$\scriptstyle{#1}$}}$}}
\def\CHAR{\mathop{\operatoratfont char \kern1pt}\nolimits}
\def\oX{{\overline{X}}}
\def\oY{{\overline{Y}}}
\def\shA{{\cal{A}}}
\def\shC{{\cal{C}}}
\def\shD{{\cal{D}}}
\def\shF{{\cal{F}}}
\def\circno#1{{\bf[#1]}}
\newcounter{alphactr}

\def\mp[[#1||#2||#3]]{\snap\hbox{$#1:\kern3pt #2\ \to\ #3$}}
\def\mapx[[#1||#2]]{\snap\hbox{$#1\ \to\ #2$}}
\def\arrowA(#1,#2)#3{\ncline[nodesep=5pt]{->}{#1}{#2}\Aput[3pt]%
    {\scriptstyle #3}}
\def\arrowB(#1,#2)#3{\ncline[nodesep=5pt]{->}{#1}{#2}\Bput[3pt]%
    {\scriptstyle #3}}
\def\arrowAB(#1,#2)#3#4{\ncline[nodesep=5pt]{->}{#1}{#2}\Aput[3pt]%
    {\scriptstyle #3}\Bput[3pt]{\scriptstyle #4}}
\def\diagramx#1{{$$\begindiagram\matrixx{#1}$$}}
\def\diagramno#1#2{{$$\begindiagram\matrixx{#2}\eqno#1$$}}
\def\begindiagram{\def\normalbaselines{\baselineskip20pt
      \lineskip3pt \lineskiplimit3pt}}
\newcounter{romanctr}
\newenvironment{romanlist}{\begin{list}{(\roman{romanctr})}{\usecounter{romanctr}}}{\end{list}} %\romanlist
\def\PLUS{+}
\def\TIMES{*}
\def\o#1{\if#1\PLUS\oplus\else\if#1\TIMES\otimes\fi\fi}

\catcode`\@=11
\def\matrixx#1{\null\,\vcenter{\normalbaselines\m@th
    \ialign{\hfil$\displaystyle ##$\hfil&&\quad\hfil$\displaystyle
##$\hfil\crcr
    \mathstrut\crcr\noalign{\kern-\baselineskip}
    #1\crcr\mathstrut\crcr\noalign{\kern-\baselineskip}}}\,}
\catcode`\@=12
\def\snap{\hskip 0pt plus 4cm\penalty1000\hskip 0pt plus -4cm}

\def\makenull#1{{\setbox0=\hbox{#1}\wd0=0pt\box0}}
\def\nulldot{\makenull{.}}

\def\mapNE#1{\hidewidth {\hphantom{\hbox{$\scriptstyle{#1}$}}}
    {\vcenter{\vbox{\hbox{\vphantom a}\hbox{$\NEarrow\kern-.3em\raise.04em\hbox
    {$\scriptstyle{#1}$}$}}}}\hidewidth}

\def\mapE#1{{\smash{\mathop{\longrightarrow}\limits^{#1}}}}
\def\mapS#1{\Big\downarrow\rlap{$\vcenter{\hbox{$\scriptstyle{#1}$}}$}}
\def\setof#1{\xmode{\{#1\}}}

\def\xmode#1{\relax\ifmmode{#1}\else{$#1$}\fi}
\catcode`\@=11
\def\operatoratfont{\operator@font}
\catcode`\@=12
\def\op{\operatoratfont}

\def\Z{\Bbb Z}
\def\O{\cal O}

\def\F{\Bbb F}
\def \I{\cal I}
\def \J{\cal J}
\def\L{\cal L}
\def\M{\cal M}
\def \N{\cal N}
\def\Q{\Bbb Q}
\def\R{\Bbb R}

\def\C{\Bbb C}
\def\P{\Bbb P}
\def\G{{\Bbb G}_{\op m}}

\def\a{{\bold a}}

\def\o*{\otimes}

\def \red{_{\op{red}}}
\def\nor{_{\op{nor}}}
\def\reg{_{\op{reg}}}

\def\mun{\mu_n}
\def\eps{\epsilon}

\def\Sing{\operatorname{Sing}}
\def\GL{\operatorname{GL}}
\def\CH{\operatorname{CH}}
\def\Gr{\operatorname{Gr}}
\def\Spec{\operatorname{Spec}}
\def\Pic{\operatorname{Pic}}
\def\Ker{\operatorname{Ker}}
\def\Cl{\operatorname{Cl}}
\def\Coker{\operatorname{Coker}}
\def\Proj{\operatorname{Proj}}

\def\dim{\operatorname{dim}}

\def\Im{\operatorname{Im}}
\def\Hom{\operatorname{Hom}}

\def\Out{\operatorname{Out}}
\def\Ext{\operatorname{Ext}}
\def\Aut{\operatorname{Aut}}
\def\Gal{\operatorname{Gal}}

\def\depth{\operatorname{depth}}
\def\a{^{\text a}}

\def\H{\operatorname{H}}
\def\D{\operatorname{D}}

\def\cat#1{\ll\negthinspace\hbox{#1}\negthinspace\gg}
\def\smallcat#1{\cat{\kern2pt\fontsize{8}{10pt}\fontshape{n}\selectfont #1\kern2pt}}

\def\DD{\mathop{\Delta}}
\def\d{\mathop{\delta}}
\begin{document}
\title[On the Picard group: torsion and the kernel]%
{On the Picard group: torsion and the kernel induced by a faithfully flat map}

\author[R. Guralnick, D. B. Jaffe, W. Raskind, and R. Wiegand]{Robert
Guralnick, David B.\ Jaffe, Wayne Raskind, and Roger Wiegand}
\address{Department of Mathematics, University of Southern California,
Los Angeles, CA 90089-1113}
\email{guralnic@mtha.usc.edu}

\address{Department of Mathematics, University of Nebraska,
Lincoln, NE 68588-0323}
\email{jaffe@cpthree.unl.edu}

\address{Department of Mathematics, University of Southern California,
Los Angeles, CA 90089-1113}
\email{raskind@mtha.usc.edu}

\address{Department of Mathematics, University of Nebraska,
Lincoln, NE 68588-0323}
\email{rwiegand@unl.edu}

\date{October 24, 1994}

\thanks{All of the authors were partially supported by
grants from the National Science Foundation.  Raskind also received support
from the NSA and from the Alexander von Humboldt-Stiftung.  Raskind thanks
A.\ Vistoli for helpful discussions.}



\maketitle

\vspace{0.1in}

\par\noindent{\bf Statement of results}

\vspace{0.1in}

For a homomorphism \mp[[ f || A || B ]] of commutative rings, let $D(A,B)$
denote the kernel of the map $\Pic(A) \to \Pic(B)$.  Let $k$ be a field and
assume that $A$ is a finitely generated $k$-algebra.

We prove a number of finiteness results for $D(A,B)$.  Here are four of them.
\circno1: Suppose $B$ is a finitely generated and faithfully flat
$A$-algebra which is geometrically integral over $k$.
If $k$ is perfect, we find that $D(A,B)$ is finitely generated.
(In positive characteristic, we need resolution of singularities to prove
this.)
For an arbitrary field $k$ of positive characteristic $p$, we find that modulo
$p$-power torsion, $D(A,B)$ is
finitely generated.  \circno2: Suppose $B = A \o*_k k^{\op sep}$.  We find
that $D(A,B)$ is finite.  \circno3: Suppose $B = A \o*_k L$, where
$L$ is a finite, purely inseparable extension.  We give examples to show
that $D(A,B)$ may be infinite.  \circno4:  Assuming resolution of
singularities,
we show that if $K/k$ is any algebraic extension, there is a finite extension
$E/k$ contained in $K/k$ such that $D(A\otimes_kE,A\otimes_kK)$ is trivial.

The remaining results are absolute finiteness results for $\Pic(A)$.
\circno5: For every $n$
prime to the characteristic of $k$, $\Pic(A)$ has only finitely
many elements of order $n$.  \circno6: Structure theorems are given for
$\Pic(A)$, in the case where $k$ is absolutely \FG.

In the body of the paper, all of these results are stated in a more general
form, valid for schemes.

\vspace{0.1in}
\par\noindent{\bf Notation and conventions}
\vspace{0.1in}

Unless otherwise stated, all rings in this paper are commutative.
Given a ring $A$, we denote by $A^*$ the group of units of $A$ and
by $A\red$ the reduction of $A$ modulo its nilradical.

Given a set $A$ and a
group $G$ acting on $A$, we let $A^G$ denote the subset of $A$ consisting of
elements fixed by $G$.  Given an abelian group $G$ and a positive integer $n$,
we let $_nG$ denote the subgroup consisting of elements whose order divides
$n$.

If $X$ is a scheme, we let $\Gamma(X)$ denote $\Gamma(X,\O_X)$.  If $X$ is a
$k$-scheme, where $k$ is a field, and $L$ is a $k$-algebra, we let $X_L$ denote
$X \times_k L$.  If $k$ is a field, $k^{\op a}$ denotes an algebraic closure of
$k$, and $k^{\op sep}$ denotes a separable closure of $k$.  We write
$X^{\op a}$, $X^{\op sep}$ instead of $X_{k^{\op a}}$, $X_{k^{\op sep}}$.

An {\it algebraic scheme\/} is a scheme which is of finite type over a field.

\section{The kernel under a separable extension}

We start by recalling some material on Galois actions, leading up to an
application of the Hochschild-Serre spectral sequence to the computation of
$\Pic$.  This material is more or less standard, but it is not available in the
literature in quite the form we need.  In particular, we want the statement of
\pref{the-kernel} to be free of noetherian hypotheses.  Later (see \ref{belch})
this will be important, because a noetherian scheme $X$ may have $\Gamma(X)$
non-noetherian, and our proofs about $\Pic(X)$ depend on understanding
$\Pic \Gamma(X)$.

If $Y$ is a scheme, and $G$ is a group (or just a set), we let $Y \times G$
denote the scheme which is a disjoint union of copies of $Y$, one for each
$g \in G$.

\begin{definition}\label{Galois-defn}
Let $X$ be a scheme, and let $Y$ be a finite \'etale $X$-scheme.  Suppose a
finite group $G$ acts on the right of $Y$ as an $X$-scheme.\footnote{Hereafter
we say simply that {\it $G$ acts on $Y/X$.}}
Then this action is {\it Galois\/} if the map
$Y \times G \to Y \times_X Y$ given by $(y,g) \mapsto (y, yg)$ is an
isomorphism of schemes.%
\footnote{If $y \in Y$, then $yg \in Y$, and $g$ induces an isomorphism of
$k(y)$ with $k(yg)$.  Therefore we get maps
\mp[[ \sigma_1, \sigma_2 || \Spec(k(y)) || Y ]], such that
$\pi \circ \sigma_1 = \pi \circ \sigma_2$, where \mp[[ \pi || Y || X ]] is the
structure map.
By the universal property of the fiber product, we obtain a morphism
\mapx[[ \Spec(k(y)) || Y \times_X Y ]], whose image is by definition the point
$(y,yg)$.}
\end{definition}

It is important to note that $G$ cannot be recovered from $Y \rightarrow X$ in
general: for instance this is the case if $Y$ consists of a disjoint union of
copies of $X$.

If we have a Galois action of $G$ on $Y/X$, and $X' \to X$ is any morphism, we
get a Galois action of $G$ on $Y \times_X X'$ as an $X'$-scheme.

Let $A \subset B$ be rings, and let a finite group $G$ act on
$B/A$, meaning that $G$ acts on (the left of) $B$ as an $A$-algebra.  Then we
shall call this action {\it Galois\/} if the action of $G$ on
$\Spec(B)/\Spec(A)$ is Galois with respect to definition \pref{Galois-defn}.
For some
definitions, stated directly for rings, see \cite{KO2, Chapter II, \S5}.
In the case where $A$ and $B$
are fields, the action of $G$ on $B/A$ is Galois if and only if $B/A$ is a
Galois extension (in the usual sense), and $G = \Aut_A(B)$.

We now recall the {\it Hochschild-Serre spectral sequence}.  This may be
found in Milne \cite{Mi1, p.\ 105}.  Although Milne refers only to locally
noetherian schemes, the argument he presents is valid for any scheme.  For
clarity, however, we note that $X\et$ as used here means the (small) \'etale
site (on an arbitrary scheme $X$), as defined in \cite{G2}.

Let $X$ be a scheme, let $Y$ be a finite
\'etale $X$-scheme, and let a Galois action of a finite group $G$ on $Y/X$
be given.   Let $\shF$ be a sheaf (of abelian groups) for the \'etale topology
on $X$.  Then the Hochschild-Serre spectral sequence is:
$$E^{p,q}_2 = H^p(G, H^q(Y\et, \pi^*\shF))
  \ \Longrightarrow \ H^{p+q}(X\et, \shF),$$%
where $\pi: Y \to X$ is the structure morphism.

Apply this with $\shF = \G$.  One has an exact sequence:
$$0 \to E^{1,0}_2 \to H^1(X\et,\G) \to E^{0,1}_2.$$%
The first term is $H^1(G,\Gamma(Y)^*)$.  The middle term is $\Pic(X)$.
The last term is $H^0(G,\Pic(Y))$, which embeds in $\Pic(Y)$.  Hence:

\begin{proposition}\label{the-kernel}
Let \mp[[ f || Y || X ]] be a finite \'etale morphism of schemes, and suppose
we have a Galois action of a finite group $G$ on $Y/X$.  Then
$$\Ker[\Pic(f)]\ \cong\ H^1(G,\Gamma(Y)^*).$$
\end{proposition}

In particular \cite{Sw2, (4.2)}, if $A \subset B$ are rings, and we are given a
Galois action of a finite group $G$ on $B/A$, then $\D(A,B) \cong \H^1(G,B^*)$.

The following example shows that one has to be a bit careful about
generalizing the proposition to the non-Galois case.

\begin{example}
Let $A$ be a complete discrete valuation ring with fraction field $K$, and let
$L$ be a finite Galois extension of $K$ with Galois group $G$.  Denote the
integral closure of $A$ in $L$ by $B$.  We would like there to be an exact
sequence:
$$0\to H^1(G,B^*)\to \Pic(A)\to \Pic(B),$$
and hence conclude that $H^1(G,B^*) = 0$.  But this is not the case in general:
$H^1(G,B^*)$ has order equal to the ramification index of $B/A$.
\end{example}

We recall some easy facts about cohomology of groups:

\begin{proposition}\label{group-cohomology-basics}
Let $G$ be a finite group and $M$ a left $\Z G$-module.
\begin{enumerate}
\item If $M$ is \FG, then $\H^n(G,M)$ is finite for every $n > 0$.
\item If $M$ has trivial $G$-action then
$H^1(G,M) \iso \Hom_{\smallcat{groups}}(G,M)$.
\end{enumerate}
\end{proposition}

\begin{proof}
(1) Since $\H^n(G,M) = \Ext^n_{\Z G}(\Z,M)$ and $\Z G$ is left Noetherian it is
clear that all the cohomology groups are finitely
generated modules (over $\Z G$ or, equivalently, over $\Z$).
Moreover, they are annihilated by $|G|$ by \cite{Br, Chap.\ III, (10.2)}.
(2) is clear from the representation of $\H^1(G,M)$ in terms of crossed
homomorphisms (or, see \cite{Bab, \S23}).  \qed
\end{proof}

\begin{proposition}\label{90}
Let a finite group $G$ act on a field $K$.  Then $H^1(G,K^*)$ is finite,
and (Hilbert's Theorem 90) it is $0$ if the group action is faithful.
\end{proposition}

\begin{proof}
Let $H$ be the subgroup of $G$ which acts trivially on $K$, and let
$\oG = G/H$.  Let $k = K^{\oG}$.  Then the extension $K/k$ is Galois, with
Galois group $\oG$.  By \pref{the-kernel},
$H^1(\oG,K^*) \iso \Ker[\Pic(k) \to \Pic(K)]$, which is $0$.  We are done
if $G$ acts faithfully on $K$.

We have the inflation-restriction exact sequence
$$0 \to \H^1(\oG, K^*) \to \H^1(G,K^*) \to \H^1(H,K^*),$$
so it is enough to show that $\H^1(H,K^*)$ is finite, which is clear from
\pref{group-cohomology-basics}(2).  \qed
\end{proof}

The next result is a variant of a well-known result due to Roquette \cite{Ro}:

\begin{proposition}\label{roquette}
Let $K$ be a field, $X$ a $K$-scheme of finite type, and  $\Lambda$ the
integral closure of $K$ in $A:=\Gamma(X)\red$.  Then $\Lambda$ is
finite-dimensional as a $K$-vector space,
and $A^*/\Lambda^*$ is a \FG\ free abelian group.
\end{proposition}

\begin{proof}
Let $\{U_1,\dots,U_m\}$ be an affine open cover of $X$, and set
$R_i = \Gamma(U_i)\red$.  Each $R_i$ is a a reduced $K$-algebra of
finite type.  We have an embedding  $A \to B := R_1\times\dots \times R_m$.
Let $P_1,\dots,P_n$ be the minimal prime ideals of $B$,
and let $C_j$ be the normalization of the domain $B/P_j$.
Then each $C_j$ is a normal domain of finite type over $K$, and
$A \subset C_1\times\dots\times C_n$.
Let $\Delta_j$ be the integral closure of $K$ in $C_j$.  By the
usual formulation of Roquette's theorem (see \cite{L, Chapter 2, (7.3)} or
\cite{Kr, (1.4)}) $C_j^*/\Delta_j^*$ is \FG\ for each $j$.  We have
$\Lambda = A \cap (\prod_j\Delta_j)$.
Therefore $A^*/\Lambda^*$ embeds in the finitely generated
group $\prod_j(C_j^*/\Delta_j^*)$.  Obviously $A^*/\Lambda^*$ is torsion-free,
and since it is finitely generated, it is free.

The fraction field $K_j$ of $C_j$ is a \FG\ field extension of $K$.
Each $\Delta_j$, being $0$-dimensional, reduced and
connected, is a field algebraic over $K$.  Since $K_j/K$
is finitely generated, so is $\Delta_j/K$.   Therefore
$\prod_j\Delta_j$ is a finite-dimensional $K$-algebra,
and hence so is its subalgebra $\Lambda$.  \qed
\end{proof}

\begin{theorem}\label{kernel2b}
Let $X$ be a scheme of finite type over a field $k$.  Let $f: Y \to X$ be a
finite, \'etale, surjective morphism of schemes.  If $k$ has positive
characteristic, assume that $X$ is reduced.  Then $\Ker[\Pic(f)]$ is finite.
\end{theorem}

\begin{remark}
If $X$ is affine, $\Pic(X) = \Pic(X\red)$, and so the assumption about
$X$ being reduced (in positive characteristic) is not needed.  In general
it is: see \pref{finite-etale-infinite}.
\end{remark}

\begin{proofnodot}
(of \ref{kernel2b}.)
We may assume that $X$ is connected.  Let $Y_0$ be a connected component
of $Y$.  Then $f|_{Y_0}$ is finite and \etale, and since $X$ is
connected, it is surjective.  Therefore we may reduce to the case where $Y$
is connected.  By \cite{Mur, (4.4.1.8)}, there
exists a scheme $Y'$ over $Y$ such that $Y' \to X$ is finite \etale\
surjective and the
action of the finite group $\Aut(Y'/X)$ on $Y'/X$ is Galois.  Therefore
\WMAT\ in fact there is a Galois action of a finite group $G$ on $Y/X$.

By \pref{the-kernel}, $\Ker[\Pic(f)] \iso H^1(G,\Gamma(Y)^*)$.  We will
complete the proof by showing that $H^1(G,\Gamma(Y)^*)$ is finite.  This
will depend only on the fact that we have a finite group $G$ acting on
an algebraic scheme $Y$, which is reduced if the characteristic is positive.

Let $B = \Gamma(Y)$.  Let $\Lambda$ be the integral closure of $k$ in $B\red$.
By \pref{roquette}, $\Lambda$ is a finite-dimensional $k$-algebra and
$B\red^*/\Lambda^*$ is finitely generated.  By
\pref{group-cohomology-basics}(1) we know that
$H^1(G,B\red^*/\Lambda^*)$ is finite.

We show that $H^1(G,\Lambda^*)$ is finite.
Write $\Lambda = \prod_{j\in J}F_j$, where each $F_j$ is a finite-dimensional
field extension of $k$ (and $J$ is a finite index set).  Since the action
of $G$ preserves idempotents, we can define an action of $G$ on $J$
by the rule $F_{gj} = gF_j$.  Let $I$ be any orbit of $G$ on $J$, and
look at $\Upsilon := \prod_{i\in I}F_i$.  Since $\Lambda^*$ is the
direct sum of the groups $\Upsilon^*$ (over the various orbits of
$G$ on $J$), it is enough to show that $\H^1(G,\Upsilon^*)$ is finite.
Fix $i \in I$, let $H$ be the isotropy subgroup of $i$, and put $F =
F_i$.  It follows from \cite{Br, Chap. III, (5.3), (5.9), (6.2)} that
$\H^1(G,\Upsilon^*) \cong \H^1(H,F^*)$, which is finite by \pref{90}.
Hence $H^1(G,\Lambda^*)$ is finite.

Running the long exact sequence of
cohomology coming from the exact sequence
$$1 \to \Lambda^* \to B\red^* \to B\red^*/\Lambda^* \to 1,$$
we conclude that $H^1(G,B\red^*)$ is finite.  Thus we are done if $X$ is
reduced.

We may assume that $\CHAR(k) = 0$.  Let $N$ be a $G$-stable nilpotent ideal
of $B$.  (For instance, we might take $N$ to be the nilradical of $B$.)  It is
enough to show that the canonical map \mapx[[ H^1(G,B^*) || H^1(G,(B/N)^*) ]]
is injective.  If $N^r = 0$, note that we can factor the map
\mapx[[ B || B/N ]] as
\diagramx{B & \mapE{} & B/N^{r-1} & \mapE{} & B/N^{r-2} & \mapE{} &
          \cdots & \mapE{} & B/N,}%
and $G$ acts on everything in the sequence, so we can reduce to the case where
$N$ has square zero.  We have an exact sequence
\diagramx{0&\mapE{} & N & \mapE{} & B^* & \mapE{} & (B/N)^* & \mapE{} & 1,}%
and it is enough to show that $H^1(G,N) = 0$.  On the one hand, $H^1(G,N)$ is
annihilated by $|G|$, and so is torsion.  On the other hand,
$H^1(G,N)$ is an $B[G]$-module, thus a $\Q$-vector space, and so is
torsion-free.  Hence $H^1(G,N) = 0$.  \qed
\end{proofnodot}

\begin{remark}\label{belch}
In the theorem, if $Y = X_L$ for some finite separable field
extension $L/k$, we will show that
there is no need to assume $X$ is reduced, even in positive
characteristic.  Indeed in that case, $\Gamma(X_L) = \Gamma(X)_L$, so we have
a Galois action of $G$ on $\Spec \Gamma(Y) / \Spec \Gamma(X)$.  Applying
\pref{the-kernel} twice, we conclude that
$\Ker[\Pic(f)] = \Ker[\Pic(\Gamma(f))]$.  But $\Pic$ of a ring is the same
as $\Pic$ of its reduction, so if $C = \Gamma(X)$, we conclude
that $\Ker[\Pic(f)] = \Ker[\Pic(C\red) \to \Pic((C\red)_L)]$.  Applying
\pref{the-kernel} once again, we see that
$\Ker[\Pic(f)] = H^1(G, (C\red)_L) = H^1(G,B\red)$.
Now the remainder of the proof of the theorem goes through.
\end{remark}

\begin{problem}\label{fidofido}
The proof of the theorem shows that if $k$ is a field, $A$ is a
\FG\ $k$-algebra (reduced if $\CHAR(k)\not=0$), and a finite group $G$ acts on
$A$, then $H^1(G,A^*)$ is finite.  Under the same hypotheses, for which
$n \in \N$ is the set $H^1(G,\GL_n(A))$ finite?
\end{problem}

\begin{theorem}\label{makes-injective}
Let $k$ be a field, and let $K/k$ be a Galois extension of fields, not
necessarily finite.  Let $S$ be a $k$-scheme of finite type, and assume that
each connected component of $S$ is geometrically connected.
Let $\Lambda$ be the integral closure of $K$ in $\Gamma(S_K)\red$.
Assume that $\Gamma(S)\red^* \Lambda^* = \Gamma(S_K)\red^*$.  Then the
canonical map $\Pic(S) \to \Pic(S_K)$ is injective.
\end{theorem}

\begin{proof}
It is enough to show that for each finite Galois extension
$L/k$ with $k \subset L \subset K$, the map $\Pic(S) \to \Pic(S_L)$ is
injective.  Let $G = \Gal(L/k)$.  Let $B_0 = \Gamma(S_L)\red$.
Let $\Lambda_0$ be the integral closure of $L$ in $B_0$.

The condition of the theorem implies that
$\Gamma(S)\red^* \Lambda_0^* = B_0^*$.  Therefore the action of $G$ on
$B_0^*/\Lambda_0^*$ is trivial.  Hence
$H^1(G,B_0^*/\Lambda_0^*) \cong \Hom_{\smallcat{groups}}(G, B_0^*/\Lambda_0^*)$
by \pref{group-cohomology-basics}(2).  Since $G$ is finite and
$B_0^*/\Lambda_0^*$ is torsion-free [by \pref{roquette}],
$H^1(G,B_0^*/\Lambda_0^*) = 0$.

Since $S$ is geometrically connected, $S_L$ is connected, and so $\Lambda_0$ is
a field.  Since $G$ acts faithfully on $\Lambda_0$,
$H^1(G,\Lambda_0^*) = 0$ by \pref{90}.

Hence $H^1(G,B_0^*) = 0$.  Arguing as in \pref{belch}, one obtains the
theorem.  \qed
\end{proof}

\begin{theorem}\label{separable-extension}
Let $k$ be a field and $S$ a scheme of finite type
over $k$.  Let $K/k$ be a separable algebraic field extension.
Then the kernel of the map $\Pic(S) \to \Pic(S_K)$ (induced by the projection
$\pi: S_K \to S$) is a finite group.
\end{theorem}

\vspace{0.08in}
\par\noindent{\bf Affine version of Theorem \ref{separable-extension}.}
\ {\it Let $k$ be a field, and let $A$ be a \FG\ $k$-algebra.  Let $K/k$ be a
separable algebraic field extension.  Then the kernel of the induced map
$\Pic(A) \to \Pic(A_K)$ is finite.}
\vspace{0.1in}

\begin{proofnodot}
(of \ref{separable-extension}).  By passing to the normal closure, we may
assume that $K/k$ is Galois (possibly infinite).  Let $B = \Gamma(S_K)\red$,
and let $\Lambda$ be the integral closure of $K$ in $B$.  By \pref{roquette},
$B^*/\Lambda^*$ is finitely generated.  Therefore we can find a finite Galois
extension $L$ of $k$ (with $L \subset K$) such that the canonical map
$(\Gamma(S_L)\red)^* \to B^*/\Lambda^*$ is surjective.  We can also choose
$L$ so that each connected component of $S_L$ is geometrically connected.  By
\pref{kernel2b} and \pref{belch}, the kernel of $\Pic(S) \to \Pic(S_L)$ is
finite, and by \pref{makes-injective}, the map
$\Pic(S_L) \to \Pic(S_K)$ is injective, so we are done.  \qed
\end{proofnodot}

We will see in \S\ref{examples-section} that \pref{separable-extension} can
fail if $K/k$ is not
assumed to be separable.  On the other hand, assuming resolution of
singularities,  we will show in \S4 that for any separated scheme $S$ of finite
type over $k$, and any algebraic field extension $K/k$,
there is an intermediate field $E/k$ of finite
degree over $k$ such that $\Pic(S_E) \to \Pic(S_K)$ is one-to-one.

We conclude this section by showing that the kernel is not always
trivial.  In fact, any finite abelian group can occur.  If, in the
following construction, one takes $K/k = \C/\R$, one
obtains the familiar example $A = \R[X,Y]/(X^2+Y^2-1)$.  The
Picard group of $A$ has order two (generated by the M\"obius
band), whereas $A\otimes_{\R}\C \cong \C[U,U^{-1}]$, which
has trivial Picard group.

\begin{example}
Let $K/k$ be a finite Galois extension with Galois
group $G$.  There is a domain $A$ of finite type over $k$ such that
$\Pic(A) \cong G/[G,G]$ but $\Pic(A\otimes_kK)$ is trivial.
\end{example}

\begin{proof}
Make $\Z G$ into a $G$-module via the left regular action.
We have an exact sequence of $G$-modules
$$0 \to \Z \mapE{\sigma} \Z G \to L \to 0,$$%
where $G$ acts trivially on $\Z$, $\sigma$ takes $1$ to $\sum_{g\in G}g$,
and $L$ is defined by the sequence.  Let $F = \Im(\sigma)$, which is the
set of fixed points of $\Z G$ under the action of $G$.  Clearly $L$ is a
free $\Z$-module of rank $|G| - 1$.

Form the group ring $B = K[L]$.  Since $L$ is a free abelian
group, $B$ is isomorphic to the Laurent polynomial ring in
$|G|-1$ variables and hence has trivial Picard group.  Note that
$G$ acts on $B$ by acting as the Galois group on $K$ and by
the $G$-module structure on $L$.  Let $A = B^G$.  Then
$A\otimes_kK \cong B$ by \cite{Sw1, (2.5)}, and the action of $G$ on
$B/A$ is Galois.  Using (\ref{the-kernel}) we get
$\Pic(A) = \D(A,B) = \H^1(G,B^*)$.  But $B^* \cong
K^* \oplus L$, so $\Pic(A) \cong \H^1(G,K^*) \oplus
\H^1(G,L) \cong \H^1(G,L)$, since $G$ acts faithfully on $K^*$ and so
$H^1(G,K^*) = 0$ by \pref{90}.  We will show that $\H^1(G,L) \cong G/[G,G]$.

We know \cite{Br, Chap. III. (6.6)} that $\H^i(G,\Z G) = 0$ for $i>0$.
Therefore, by the long exact sequence of cohomology we have
$\H^1(G,L) \cong \H^2(G,\Z)$.  But
$\H^2(G,\Z) \cong \H^1(G,\Q/\Z) \cong
 \Hom_{\smallcat{groups}}(G, \Q/\Z) \cong G/[G,G]$.
(See \cite{Bab, \S23}.)  \qed
\end{proof}

\section{Torsion in Picard groups}\label{torsion-section}

We note that $_n\Pic(R)$ can be infinite even for $R$ an domain \FG\ over
a field.  For example, take $R = k[T^2,T^3]$, where $k$ is an
infinite field of  characteristic $p > 0$.  Then $\Pic(R)$ is
isomorphic to to the additive group of $k$ and is therefore an
infinite group of exponent $p$.   As long as we avoid the
characteristic, however, this cannot happen.  First we need the following
lemma (cf.\ \cite{Bas1, IX, (4.7)}):

\begin{lemma}\label{kernel-is-torsion}
Let \mp[[ f || X || S ]] be a finite flat morphism of schemes, of constant
degree $d > 0$.  Then $\Ker[\Pic(f)]$ is $d$-torsion.
\end{lemma}

\begin{proof}
Let $\M \in \Pic(S)$.  Then $f_*f^*\M \iso \M \o* f_* \O_X$ as
$\O_S$-modules.  If moreover $\M \in \Ker[\Pic(f)]$, then
$f_* f^* \M \iso f_* \O_X$.  Hence $\M \o* f_* \O_X \iso f_* \O_X$.  Apply
$\wedge^d$, yielding
$\M^{\o* d} \o* \wedge^d(f_*\O_X) \iso \wedge^d(f_*\O_X)$, and hence
$\M^{\o* d} \iso \O_S$.  \qed
\end{proof}

Now let us generalize to the proper case.  First we need:

\begin{lemma}\label{is-locally-free}
Let \mp[[ f || X || S ]] be a proper flat morphism of noetherian
schemes.  Then $f_*\O_X$ is a locally free $\O_S$-module.
\end{lemma}

\begin{proof}
We may assume that $S$ is affine.  Let $W = \SPEC(f_*\O_X)$
(cf.\ \cite{Ha, II, exercise 5.17}), and label morphisms
\diagramx{X & \mapE{\varphi} & W & \mapE{h} & S.}%
Since $f$ is proper, $f_*\O_X$ is coherent \cite{EGA$3_1$, 3.2.1}, and so it is
enough to show that
$h$ is flat.  Therefore it is enough to show that for any injection $i$ of
coherent $\O_S$-modules, $h^*(i)$ is also injective.  By
construction, $\varphi_* \O_X = \O_W$, and so
$\varphi_* \varphi^* h^*(i) = h^*(i)$ by the projection formula
\cite{Ha, II, exercise 5.1d}.  Since $f$ is
flat, $f^*(i)$ [which equals $\varphi^* h^*(i)$] is injective.  Since
$\varphi_*$ is left exact, $\varphi_* \varphi^* h^*(i)$ is also injective.
\qed
\end{proof}

\begin{corollary}\label{proper-bounded}
Let \mp[[ f || X || S ]] be a surjective proper flat morphism of noetherian
schemes.  Then $\Ker[\Pic(f)]$ is a bounded torsion group.
\end{corollary}

\begin{proof}
We may assume that $S$ is connected.
By \pref{is-locally-free}, $f_*\O_X$ is a locally free $\O_S$-module.
Factor $f$ as in the proof of \pref{is-locally-free}.  By the projection
formula, $\varphi_* \varphi^*$ is the identity, so $\Pic(\varphi)$ is
injective.
Apply \pref{kernel-is-torsion}.  \qed
\end{proof}

\begin{theorem}\label{finitely-many}
Let $S$ be a scheme of finite type over a field $k$ and let $n \in \NN$ be
prime to the characteristic of $k$.  Then $_n\Pic(S)$ is finite.
\end{theorem}

\vspace{0.03in}
\par\noindent{\bf Affine version of Theorem \ref{finitely-many}.}
\ {\it Let $k$ be a field, and let $A$ be a \FG\ $k$-algebra.  Let $n \in \NN$
be prime to the characteristic of $k$.  Then $_n\Pic(A)$ is finite.}
\vspace{0.1in}

\begin{proofnodot}
(of \ref{finitely-many}).  Suppose first that $k$ is separably closed.
The argument in this case seems to be fairly well known and was pointed
out to us several years ago by David Saltman and Tim Ford.
We consider the Kummer sequence \cite{SGA4$1\over2$, p.\ 21, (2.5)}:
\diagramx{1 & \mapE{} & \mun & \mapE{} & \G & \mapE{n} & \G & \mapE{} & 1}%
This is an exact sequence of sheaves for the \'etale topology
on $S$.  Taking \'etale cohomology, we get an exact sequence
\diagramno{(*)}{\H^1(S\et,\mun) & \mapE{} & \H^1(S\et,\G) & \mapE{n} &
                \H^1(S\et,\G).}%
Now $\H^1(S\et,\G)\cong\Pic(S)$ by \cite{SGA4$1\over2$, p.\ 20, (2.3)}.
Since $k$ is separably closed, $\mun$ is isomorphic to the
constant sheaf $\Z/n\Z$.  (See \cite{SGA4$1\over2$, p.\ 21, (2.4)}.)
Therefore $\mun$ is constructible \cite{SGA4$1\over2$, p.\ 43, (3.2)}.  By
\cite{SGA4$1\over2$, p.\ 236, (1.10)} $\H^1(S\et,\mun)$ is finitely generated.
Hence the kernel of the map $n$ in $(*)$ is finite, so the theorem holds when
$k$ is separably closed.  Apply \pref{separable-extension}.  \qed
\end{proofnodot}

As a consequence of \pref{proper-bounded} and \pref{finitely-many}, we have:

\begin{corollary}\label{proper-Pic-finite}
Let $S$ be a scheme of finite type over a field $k$.  Let \mp[[ f || X || S ]]
be a proper flat surjective morphism of schemes.  Assume that over each
connected component of $S$, the rank of the locally free sheaf $f_*\O_X$
is invertible in $k$.  Then $\Ker[\Pic(f)]$ is finite.
\end{corollary}

We can use (\ref{finitely-many}) to answer a question of S.\ Montgomery about
outer automorphism groups of Azumaya algebras:

\begin{corollary}
Let $k$ be a field.  Let $R$ be a \FG\ $k$-algebra.  Let $A$ be a (not
necessarily commutative) Azumaya algebra over $R$ of degree $d$.  Then
$\Out_R(A)$ is finite whenever $d$ is not a multiple of
the characteristic of $k$.
\end{corollary}

\begin{proof}
There is an embedding (see e.g.\ \cite{DI}) of $\Out_R(A)$ in $\Pic(R)$,
and in fact the image of $\Out_R(A)$ is contained in $_d \Pic(R)$.
(See \cite{KO1}.)  \qed
\end{proof}

In fact, one can show that the finiteness of $_{d^e}\Pic(R)$ for all $e$
is equivalent to the finiteness of $\Out_R(M_{d^e}(R))$ for
all $e$ (see \cite{BG}).

Finally, we consider $n$-torsion in a normal algebraic scheme.  It turns
out that this is finite, even if $n$ is not prime to the characteristic,
at least assuming that resolution of singularities holds.  To prove this, we
need to know what happens to $\Pic$ of a normal
scheme when its singular locus is deleted.  The kernel is described by the
following ``folklore'' lemma, which we state in greater generality for later
application \pref{seminormal-S2}:

\begin{lemma}\label{S2}
Let $X$ be a noetherian $S_2$ scheme%
\footnote{Recall that a noetherian scheme $X$ is by definition $S_2$ if for
every $x \in X$, $\depth \O_{X,x} \geq \min\setof{\dim \O_{X,x},2}$.},
and let $C \subset X$ be a closed subset
of codimension $\geq 2$.  Let $U = X - C$.  Then the canonical map
$\Pic(X) \to \Pic(U)$ is injective.
\end{lemma}

\begin{proof}
Let $\L$ be a line bundle on $X$ which becomes trivial on $U$.  For any
line bundle $\M$ on $X$, the long exact sequence of local cohomology gives us
$$\H^0_C(X,\M) \to \H^0(X,\M) \mapE{\rho_{\M}} \H^0(U,\M) \to \H^1_C(X,\M).$$%
Since $X$ is $S_2$, the end terms vanish \cite{G1, (1.4), (3.7), (3.8)}.
Let $\phi: \O_U \to \L|_U$  be an isomorphism.  Since $\rho_{\L}$
is an isomorphism, we can lift $\phi$ to a morphism
$\psi: \O_X \to \L$.  Since $\rho_{\L^{-1}}$ is an isomorphism,
we can lift $\phi^{-1}$ to a morphism $\psi': \L \to \O_X$.  Since
$\rho_{\O_X}$ and $\rho_{\L}$ are isomorphisms, $\psi' \circ \psi$
and
$\psi \circ \psi'$ are the identity maps.  \qed
\end{proof}

\begin{lemma}\label{alg-closed-regular}
Let $k$ be an algebraically closed field.  Assume that resolutions of
singularities exist for varieties over $k$.  Let $X$ be a normal $k$-scheme
of finite type.  Let $n \in {\Bbb N}$.  Then $_n \Pic(X)$ is finite.
\end{lemma}

\begin{sketch}
We may assume that $X$ is connected.  By \pref{S2}, we may replace $X$ by
$X_{\op reg}$ and so assume that $X$ is regular.  If we further replace
$X$ by a nonempty open subscheme, we kill a \FG\ subgroup of $\Pic(X)$.  In
this way we may reduce to the case where $X$ is affine.  Since we have
resolution of singularities, we can embed $X$ as an open subscheme of a
regular $k$-scheme $\oX$.  Then $\Pic(X)$ is the quotient of $\Pic(\oX)$
by a \FG\ subgroup.  Now $\Pic^0(\oX)$ is the group $A(k)$ of $k$-valued points
of an abelian variety $A$ over $k$, and $\Pic(\oX)/\Pic^0(\oX)$ is
\FG\ (see e.g.\ \cite{K, (5.1)}).  Therefore it suffices to show that $_n A$ is
finite.  This is well-known \cite{Mum, p.\ 39}.  \qed
\end{sketch}

\section{Faithfully flat extensions}

The main results of this section are \pref{faithfully-flat-fg},
\pref{faithfully-flat-fg-imperfect}, and \pref{open-cover}, which give
information
about the kernel of the map on Picard groups induced by a faithfully flat
morphism of algebraic schemes.  First we consider the case where the target of
the morphism is normal, in which case we can weaken the hypothesis of faithful
flatness.

\begin{theorem}\label{normal-faithfully-flat}
Let $k$ be a field.  Let $X$ be a normal $k$-scheme of finite type.  Let
\mp[[ f || Y || X ]] be a dominant morphism of finite type.
Then $\Ker[\Pic(f)]$ is \FG\ (if $\CHAR(k) = 0$) and is the
direct sum of a \FG\ group and a bounded $p$-group (if $\CHAR(k) = p > 0$).
If $k$ is algebraically closed, and resolution of singularities holds, then
$\Ker[\Pic(f)]$ is \FG.
\end{theorem}

\begin{proof}
We may assume that $X$ is connected.  By \pref{S2}, \WMAT\ $X$ is regular.
Then we may replace $X$ by any nonempty open subscheme.  In particular,
\WMAT\ $X$ is affine.
Moreover, by replacing $Y$ by a suitable open subscheme, \WMAT\ $Y$
is affine too.  We may assume that $Y$ is a regular integral scheme.  We
may embed $Y$ as an open subscheme of an $X$-scheme $\oY$ which is projective
over $X$ and is an integral scheme as well.  Replace $\oY$ be its
normalization.
Now by again replacing $X$ by a nonempty open subscheme, we may
assume (by generic flatness) that the morphism \mapx[[ \oY || X ]] is flat;
certainly \WMAT\ it is surjective.  Call this morphism $\varphi$.  By
\pref{proper-bounded},
$\Ker[\Pic(\varphi)]$ is a bounded torsion group.  Hence by
\pref{finitely-many}, $\Ker[\Pic(\varphi)]$ is finite (if $\CHAR(k) = 0$) and
is the direct sum of a finite group and a bounded $p$-group
(if $\CHAR(k) = p > 0$).  [By \pref{alg-closed-regular}, if $k$ is
algebraically closed, then $\Ker[\Pic(\varphi)]$ is always finite.]
By \pref{S2} and \cite{Ha, II, (6.5c), (6.16)},
$\Ker[\Pic(\oY) \rightarrow \Pic(Y)]$ is \FG.  The theorem follows.  \qed
\end{proof}

Now we want to see what happens when we consider a faithfully flat morphism
\mapx[[ Y || X ]], where $X$ is not necessarily normal.  The issue is
complicated by nilpotents, even in the affine case.  The problem [see examples
\pref{reduction-not-faithfully-flat}, \pref{rnff2} below] is that one can have
a domain $A$, and a \FG\ faithfully flat $A$ algebra $B$, such that $B\red$ is
not flat over $A$.

\begin{lemma}\label{B-cap-K-if-flat}
Let $A$ be a reduced ring, with total ring of fractions $K$.  Let $B$ be a
faithfully flat $A$-algebra.  Then (inside $B \o*_A K$) $B \cap K = A$.
\end{lemma}

\begin{proof}
Let $b \in B \cap K$, so we have an equation of the form $bu = v$, for some
$u, v \in A$ with $u$ a non-zero-divisor.  By the equational criterion for
faithful flatness \cite{Bour, Ch.\ I, \S3, $\op{n}^\circ$ 7, Prop.\ 13},
$b \in A$.  \qed
\end{proof}

\begin{example}\label{reduction-not-faithfully-flat}
Let $k$ be a field of characteristic $2$.  Let $A = k[s,t]/(s^2-t^3)$.
Let $B = A[x,y]/(x^2-s,y^2-t)$.  Then $B$ is a faithfully flat $A$-algebra.
Now $B \iso k[x,y]/(x^4-y^6)$, so the nilradical of $B$ is generated by
$(x^2-y^3)$.  Hence $B\red = A[x,y]/(x^2-s,y^2-t,x^2-y^3)$.  Since
$s=ty$ in $B\red$, we have $y \in B\red \cap A\nor$ and $y \notin A$, so it
follows from \pref{B-cap-K-if-flat} that $B\red$ is not flat over $A$.
One can also prove this directly by taking $M = A/(t)$ and checking that
the map \mapx[[ M || M \o*_A B\red ]] is not injective.
\end{example}

\begin{example}[shown to us by Bill Heinzer and Sam Huckaba]\label{rnff2}
It is known that there is a smooth curve $C \IN \CP^3$ of degree $8$ and
genus $5$ which is set-theoretically the intersection of two surfaces $S$, $T$,
but which is not arithmetically Cohen-Macaulay (see \cite{Bar}, \cite{Hu}).
Let $\tC = S \cap T$,
scheme-theoretically.  Choose lines $L, L'$ which are noncoplanar and do not
meet $C$.  Then projection from $L$ onto $L'$ defines a
nonconstant morphism \mapx[[ \tC || \P^1 ]].
Let $A$ be the homogeneous coordinate ring of $\P^1$, and let
$B$ be the homogeneous coordinate ring of $\tC$.  Then $A = \C[x,y]$, $B$
is Cohen-Macaulay, $A \IN B$, and $B$ is module-finite over $A$, so $B$ is
faithfully flat over $A$.  (See \cite{Ma, p.\ 140}.)  On the other hand,
$B\red$ is not \CM, so $B\red$ is {\it not\/} flat over $A$.
\end{example}

It is not clear if the behavior illustrated by the characteristic $p$
example can be mimicked in characteristic zero:

\begin{problem}\label{furry-friend}
Let $k$ be an algebraically closed field of characteristic zero.  Let
$A$ be a \FG, reduced $k$-algebra.  Let $B$ be a faithfully flat and
\FG\ $A$-algebra.  Do we have $B\red \cap A\nor = A$?
\end{problem}

\begin{theorem}\label{faithfully-flat-fg}
Let $k$ be a perfect field.  Assume that resolutions of singularities exist
for varieties over $k^a$.  Let $X$ and $Y$ be geometrically integral
$k$-schemes of finite type.  Let \mp[[ f || Y || X ]] be a faithfully flat
morphism of $k$-schemes.  Then $\Ker[\Pic(f)]$ is \FG.
\end{theorem}

We are not sure to what extent the hypothesis ``geometrically integral''
can be relaxed.  Certainly if $X$ is not affine, nonreduced, and $Y$ is
disconnected, one can have trouble \pref{nonreduced-cover-example}.
In positive characteristic, the
assumption that $Y$ is reduced is needed \pref{random-rabbits}.
In characteristic zero, we do not know if it is necessary to assume that $Y$ is
reduced.  [If the answer to \pref{furry-friend} is yes, then we do not
need to assume $Y$ reduced.]  We do not know if it is necessary to
assume that $X$ and $Y$ are geometrically irreducible.  Cf.\ \pref{open-cover}.

\vspace{0.1in}
\par\noindent{\bf Affine version of Theorem \ref{faithfully-flat-fg}.}
\ {\it Let $k$ be a perfect field.  Assume that resolutions of singularities
exist for varieties over $k^a$.  Let $A$ be a \FG\ $k$-algebra.
Let $B$ be a \FG\ and faithfully flat $A$-algebra, which is geometrically
integral over $k$.  Then $\Ker[\Pic(A) \rightarrow \Pic(B)]$ is \FG.}
\vspace{0.1in}

\begin{proofnodot}
(of \ref{faithfully-flat-fg}).
By \pref{separable-extension}, $\Ker[\Pic(X) \rightarrow \Pic(X^{\op a})]$ is
finite, so \WMAT\ $k$ is algebraically closed.

By \pref{normal-faithfully-flat}, the theorem holds when $X$ is normal.  To
complete the proof, we need to show that
$\Ker[\Pic(f)] \cap \Ker[\Pic(X) \to \Pic(X\nor)]$ is \FG.  If we relax the
assumptions on $Y$, by assuming only that each connected component of $Y$ is
geometrically integral, it is enough to do the two cases:
\begin{romanlist}
\item $Y$ is an ``open cover'' of $X$;
\item $X$, $Y$ are both affine.
\end{romanlist}

Let $P$ be the fiber product of $X\nor$ and $Y$ over $X$. Then $P$ is reduced.
Let \mp[[ \pi || X\nor || X ]] and \mp[[ \tau || P || Y ]] be the canonical
maps.  Let $\shC$ be the quotient sheaf $(\pi_*\O_{X\nor}^*)/\O_X^*$, and
similarly let $\shD = (\tau_*\O_P^*)/\O_{Y}^*$.  (Note that the canonical map
\mapx[[ \O_{Y} || \tau_*\O_P ]] is injective.)
We have a commutative diagram with exact rows (and some maps labelled):
\diagramx{1 & \mapE{} & {\Gamma(X\nor)^*/\Gamma(X)^*} & \mapE{i}
                           & \Gamma(\shC) & \mapE{}
                           & \Ker[\Pic(X) \to \Pic(X\nor)] & \mapE{} & 1\cr
          && \mapS{\lambda} && \mapS{\delta} && \mapS{}\cr
          1 & \mapE{} & {\Gamma(P)^*/\Gamma(Y)^*} & \mapE{} & \Gamma(\shD)
                      & \mapE{}
                      & \Ker[\Pic(Y) \to \Pic(P)] & \mapE{} & 1\nulldot}%
We will be done with the proof if we can show that $\delta$ is injective
and that $\Coker(\lambda)$ is \FG.

We show that $\delta$ is injective.  In case (i) this is clear, since
$\shC$ is a sheaf.  So \WMAT\ $X$, $Y$ are both affine.
Let $A = \Gamma(X)$, $B = \Gamma(Y)$.  Let $\alpha \in \Gamma(\shC)$.
Then there exist elements $\vec f1n \in A$ with $(\vec f1n) = (1)$ and
elements $\alpha_i \in (A\nor)_{f_i}^*$ $(i=1,\ldots,n)$,
$\beta_{ij} \in A_{f_i f_j}^*$ such that
$\alpha_i = \beta_{ij}\alpha_j$ in $(A\nor)_{f_i f_j}^*$ and $\alpha$
induces $\vec \alpha1n$.  Suppose $\alpha$ maps to $1 \in \Gamma(\shD)$.
Then for each $i$ we have elements $b_i \in B_{f_i}^*$ such that
$b_i = \alpha_i$ in $B_{f_i} \o*_A A\nor$ for each $i$.  By
\pref{B-cap-K-if-flat}, applied with $A_{f_i}$ substituted for $A$, we see that
$\alpha_i \in A_{f_i}$.  Hence $\alpha = 1$.  Hence $\delta$ is injective.

To complete the proof, we need to show that
$$\Coker(\lambda) = {\Gamma(P)^* \over \Gamma(Y)^* \Gamma(X\nor)^*}$$%
is \FG.  For this, we may assume that $Y$ is connected.  By \pref{roquette}, it
is enough to show that the map \mapx[[ P || X\nor ]] is bijective on connected
components.  This is true since $X$ and $Y$ are integral schemes.  \qed
\end{proofnodot}

We note that the proof breaks down if we do not assume that $Y$
is geometrically irreducible.  Indeed, without this hypothesis, the
map \mapx[[ P || X\nor ]] may not be bijective on connected components:

\begin{example}
Let $k$ be a field of characteristic $\not= 2$.  Let $A$ be the subring
$k[t^2, t+t^{-1}]$ of
$k[t,t^{-1}]$.  Let $B = A[x]/(x^2-t^2)$.  Then $B$ is a faithfully flat
$A$-algebra.  We will see that (i) $\Spec(B)$ is connected, but (ii)
$\Spec(B \o*_A A\nor)$ is not connected.  First note that
$t = (t^2+1)/(t+t^{-1})$, from which it follows that $A\nor = k[t,t^{-1}]$.
Then $B \o*_A A\nor \iso A\nor \times A\nor$, so (ii) holds.
Let $z = t^2 + 1$, $y = t+t^{-1}$.  Then $A \iso k[y,z]/(y^2+z^2-y^2z)$, so
$$B \iso k[x,y,z]/(x^2-z+1, y^2+z^2-y^2z) \iso k[x,y]/((x^2+1)^2-x^2y^2).$$%
{}From this one sees that $\Spec(B)$ has two smooth components, meeting at
the two points $(\pm i,0)$.  Hence $\Spec(B)$ is connected.
\end{example}

\begin{theorem}\label{faithfully-flat-fg-imperfect}
Let $k$ be a field of characteristic $p > 0$.  Let $X$ and $Y$ be
geometrically irreducible $k$-schemes of finite type.  Let
\mp[[ f || Y || X ]] be a faithfully flat morphism of $k$-schemes.
Then $\Ker[\Pic(f)]$ is the direct sum of a \FG\ group and a $p$-group.
\end{theorem}

\vspace{0.05in}
\par\noindent{\bf Affine version of Theorem
\ref{faithfully-flat-fg-imperfect}.}
\ {\it Let $k$ be a field of characteristic $p > 0$.  Let $A$ be a
\FG\ $k$-algebra.  Let $B$ be a \FG\ and faithfully flat $A$-algebra such
that $\Spec(A \o*_k k^a)$ is irreducible.  Then the group
$\Ker[\Pic(A) \rightarrow \Pic(B)]$
is the direct sum of a \FG\ group and a $p$-group.}
\vspace{0.1in}

\begin{proofnodot}
(of \ref{faithfully-flat-fg-imperfect}).
It follows from \pref{kernel-is-torsion} and \pref{finitely-many} that
$\Ker[\Pic(X) \to \Pic(X^a)]$ is the direct sum of a \FG\ group and a
$p$-group, so \WMAT\ $k$ is algebraically closed.  Since
$\Ker[\Pic(X) \to \Pic(X\red)]$ is a $p$-group, \WMAT\ $X$ is reduced.

If $Y$ is reduced, we are done by \pref{faithfully-flat-fg}.  In the
general case, we modify the argument of the proof of
\pref{faithfully-flat-fg}.  Since the scheme $P$ in that argument may be
nonreduced, we do not get that $\Coker(\lambda)$ is \FG, but rather
$(*)$ it is \FG\ mod $p$-power torsion.  We still get that
$$M := {\Gamma(P)\red^* \over \Gamma(Y)\red^* \Gamma(X\nor)^*}$$%
is \FG.  Let $N$ be the nilradical of $\Gamma(P)$.  Then the kernel of the
canonical map \mapx[[ \Coker(\lambda) || M ]] is a quotient of $1+N$, and
hence is a bounded $p$-group, so $(*)$ follows.  \qed
\end{proofnodot}

\begin{remark}
By \pref{bounded}, it will follow that the $p$-group of the
theorem is actually a bounded $p$-group, provided that resolution of
singularities is valid.
\end{remark}

\begin{theorem}\label{open-cover}
Let $X$ be a reduced scheme of finite type over a perfect field.  Let
$\vec U1n$ be an open cover of $X$.  Then the canonical map
\dmapx[[ \Pic(X) || \Pic(U_1) \times \cdots \times \Pic(U_n) ]]
has \FG\ kernel.
\end{theorem}

For the case where $X$ is normal (in fact, any noetherian normal scheme), this
is easy, following roughly from \pref{S2}.
For the general case, one can follow the proof of \pref{faithfully-flat-fg},
taking $Y$ to be the disjoint union of the $U_i$'s;
one need only adjust the last sentence.

\section{Examples where the kernel is not finitely generated}%
\label{examples-section}

Let \mp[[ f || Y || X ]] be a faithfully flat morphism of noetherian schemes.
The results \ref{kernel2b}, \ref{separable-extension}, \ref{proper-Pic-finite},
\ref{normal-faithfully-flat}, \ref{faithfully-flat-fg}, and \ref{open-cover}
all give conditions under which $\Ker[\Pic(f)]$ is \FG.  While it is
reasonable to think that there are unifying and more general results with
the same conclusion, we do not know what form such results should take.
With this in mind, we give in this section a varied collection of examples
in which $\Ker[\Pic(f)]$ is {\it not\/} \FG.

We mention an obviously related question, about which we know very little.
For which faithfully flat proper morphisms \mp[[ f || Y || X ]] is it the
case that the map
\dmapx[[ \setofh{iso.\ classes of vector bundles on $X$} ||
         \setofh{iso.\ classes of vector bundles on $Y$} ]]%
is finite-to-one?  Cf.\ \pref{fidofido}, \pref{woofwoof}.

Returning to $\Pic$,
first we see that there are examples with $X$ algebraic of positive
characteristic, and $f$ a finite \etale\ morphism.  For these, by
\pref{kernel2b}, $X$ must be nonreduced.

\begin{example}\label{finite-etale-infinite}
Let $k$ be an algebraically closed field of characteristic $p > 0$.  Let
\mp[[ \lambda || F || E ]] be a finite \etale\ morphism of varieties over $k$,
of degree $p$, and suppose we have a Galois action of a finite group $G$
on $F/E$.  Let $R = k[\eps]/(\eps^p)$.
Let $X = E \times_k R$, $Y = F \times_k R$, and let \mp[[ f || Y || X ]] be
the induced morphism.  Then $\Ker[\Pic(f)] \iso H^1(G, \Gamma(Y)^*)$ by
\pref{the-kernel}.  Since $\Gamma(Y)$ is just $R$, and $G$ acts trivially on
it, we have by \pref{group-cohomology-basics}(2) that
$\Ker[\Pic(f)] \iso \Hom(\Z/p\Z, R^*)$, which is not \FG, since
$(1+c\eps)^p = 1$ in $R$ for each $c \in k$.

To get examples of such morphisms $\lambda$,
let $E$ be an elliptic curve over $k$ which is not supersingular.  Then
the Tate module $T_p(E)$ is isomorphic
to the $p$-adic integers $\Z_p$.  According to a theorem of Serre-Lang
\cite{SGA1, XI, (2.1)}, the $p$-primary part of the algebraic fundamental
group $\pi_1(E)$ is $T_p(E)$.  Therefore there exists a surjective
homomorphism \mapx[[ \pi_1(E) || \Z/p\Z ]], and so there exists a variety $F$
and a morphism $\lambda$ as indicated.  More concretely, this may be seen
as follows.  The multiplication by $p$ map $E\ \mapE{p}\ E$ factors through
$E^{(p)}$, the scheme defined in essence by raising the coefficients in
the equation defining $E$ to the \th{p} power.  The induced map
\mapx[[ E^{(p)} || E ]], called the {\it Verschiebung}, is exactly $\lambda$.
\end{example}

Now we see that there are examples with $X$ reduced and algebraic, even over an
algebraically closed field (of positive characteristic),
and $f$ a finite flat morphism:

\begin{example}\label{random-rabbits}
Take example \pref{reduction-not-faithfully-flat}, and use $X = \Spec(A)$,
$Y = \Spec(B)$.  Since the map \mapx[[ A || B\red ]] factors through $A\nor$,
it follows that $\Ker[\Pic(f)] = k$, which is not \FG\ if $k$ is infinite.
\end{example}

Now we see that there are examples with $X$ algebraic and $Y$ \etale, even
in characteristic zero.  Indeed, one may take $Y$ to be an ``affine open
cover''
of $X$:

\begin{example}\label{nonreduced-cover-example}
Let $T$ be a projective variety over an algebraically closed field $k$,
and let $\shF$ be a coherent sheaf on $T$ with $H^1(T,\shF) \not= 0$.
Make $\shA := \O_T \oplus \shF$ into an $\O_T$-algebra by forcing
$\shF \cdot \shF = 0$.  Let $X = \SPEC(\shA)$.  Let $Y$ be an affine open
cover of $X$, i.e.\ the disjoint union of the schemes in such a cover.
Then $\Ker[\Pic(X) \to \Pic(X\red)] \IN \Ker[\Pic(X) \to \Pic(Y)]$, and
$\Ker[\Pic(X) \to \Pic(X\red)]$ is a nonzero vector space (over $k$), so
$\Ker[\Pic(X) \to \Pic(Y)]$ is not \FG.
\end{example}

Now we give families of examples in which $X$ and $Y$ can be chosen to be
reduced noetherian schemes of characteristic zero, but the morphism is not of
finite type:

\begin{example}
Given any ring $A$ there exists a faithfully flat extension $B$ with
$\Pic(B) =1$: take $B=A[x]$ localized at the set of primitive polynomials
(i.e.\ polynomials such that the coefficients generate the unit ideal).
For the fact that $\Pic(B)=1$, see \cite{EG, (5.4), (3.5), (2.6 with $R=S$)}.
\end{example}

\begin{example}
Let $k$ be a field and let $X = \Spec(k[t^2,t^3])$.  Let
$$Y = \Spec(k[t,t^{-1}]) \times k[t^2,t^3]_{(t^2,t^3)}),$$%
which is the disjoint union of $X\reg$ and $\Spec \O_{X,x}$, where $x$ is
the singular point of $X$.  Then $\Ker[\Pic(f)] = k$.
\end{example}

Now we will give examples based on purely inseparable base extension.  For
purely inseparable extensions we cannot use Galois cohomology to control the
Picard group.  Instead, we use differentials, following the ideas in Samuel's
notes on unique factorization domains \cite{Sam}. We build a purely inseparable
form of the affine line whose Picard group is infinite.

We need the following version of a result of Samuel \cite{Sam, 2.1, p.\ 62}:

\begin{lemma}\label{son-of-samuel}
Let $B$ be a domain of characteristic $p > 0$ with fraction field $Q(B)$.  Let
$\d$ be a $\Z$-linear
derivation of $B$ with $A$ the subring of invariants of $\d$ (i.e.\ the
elements with $\d(a)=0$).  Let $\DD$ be the logarithmic derivative
of $\d$ defined on $Q(B)^*$ (i.e. $\DD(b)=\d(b)/b)$.  If $M_1, M_2$ are
invertible ideals of $A$ with $M_iB = b_iB$ for each $i$ ($b_1,b_2 \in B$),
then $M_1 \cong M_2$ if and only if $\Delta(b_1)-\Delta(b_2) \in \Delta(B^*)$.
\end{lemma}

\begin{proof}
If $M_1 \cong M_2$, then $b_1=aub_2$ for some $a$ in the quotient
field of $A$ and $u \in B^*$.  Thus, $\Delta(b_1)=\Delta(b_2) + \Delta(u)$.

Conversely, if $\Delta(b_1)-\Delta(b_2) = \Delta(u)$ for some
$u \in B^*$, then replacing $b_1$ by $b_1u$ allows us to assume
that $\Delta(b_1)=\Delta(b_2)$ whence $\delta(b_1^{-1}b_2 )=0$.
Set $a=b_1^{-1}b_2$.  Thus, $a$ is in the quotient field of $A$.
Replacing $M_1$ by $aM_1$ allows us to assume that $M_1B=M_2B$.
We claim that this implies $M_1=M_2$.  It suffices to check this
locally and so we may assume that each $M_i$ is principal.  Let
$a_i$ be a generator for $M_i$.  It follows that $a_1/a_2$ is a unit
in $B$.  We have $\delta(a_1/a_2) = 0$ and $a_1/a_2 \in B$, so $a_1/a_2 \in A$.
Similarly, $a_2/a_1 \in A$, so $a_1/a_2 \in A^*$, and hence $M_1 = M_2$.  \qed
\end{proof}

\begin{example}\label{deranged-derivations}
Let $k$ be a separably closed imperfect field of
characteristic $p$.  Let $\alpha \in k - k^p$.  Let
$q$ be a power of $p$ with $q > 2$, and let $K = k(\alpha^{1/q})$.
Let $A = k[X,Y]/(X^q-X-\alpha Y^q)$.  Set $r=q/p$.  Then $A$ is a
Dedekind domain and $P=\Pic(A)$ is a group of exponent $q$ with $P/{_r}P$
infinite.  Also, $A_K \iso K[Z]$.  Hence the kernel of the map
$\Pic(A) \rightarrow \Pic(A_K)$ is infinite of exponent $q$.
\end{example}

\begin{proof}
Let $\beta = \alpha^{1/q}$.  In the polynomial ring  $B = K[Z]$, put $x = Z^q$
and $y = \beta^{-1}(Z^q - Z)$.  Note that $Z = x-\beta y$, so that
$K[Z] = K[x,y]$. Moreover, we can identify $A$ with the subring $k[x,y]$ of
$B$.  Then $A_K = B$.  Thus $A$ is a Dedekind domain and $\Pic(A_K)$ is
trivial.  It follows that $\Pic(A)$ has exponent dividing $q$.

Let $V=\{a \in k : a^q -a \in \alpha k^q\}$.  Since $k$ is separably closed,
for every $b \in k$ there is an element $a\in V$ such that
$$a^q - a -\alpha b^q = 0.$$
Thus $V$ is infinite.  Given $a \in V$, define $b \in k$
by the above equation.  Let $M(a)$ be the maximal ideal $(x-a,y-b)$ of $A$.
Set $c= a- \beta b=a^{1/q}$. Note that $M(a)B=(Z-c)B$.

Let $a_1$ and $a_2$ be distinct elements of $V$.  Let $b_i$ and $c_i$ be the
corresponding elements defined above.  It suffices to  show that
$M(a_1)^r$ and  $M({a_2})^r$ are nonisomorphic.

Let $\gamma = \alpha^{1/p}$.  Define a derivation $\d$ on $k[\gamma]$ with
$\d(k)=0$ and $\d(\gamma)= \gamma$.  Extend this derivation to $A_0=A[\gamma]$
(and to its quotient field) by taking $\d$ trivial on $A$.  Set $W=Z^r$.  Since
$Z=x - \beta y$, $W=x^r -\gamma y^r$.  Thus, $\d(W)= -\gamma y^r = W-W^q$.  Let
$\DD$ denote the logarithmic derivative of $\d$.

The following observation will be useful:  $\d(c_i^r)=-\gamma b_i^r=
c_i^r -c_i^{rq}$. Thus
$$\DD(W-c_i^r) ={{W-W^q-c_i^r+c_i^{rq}}\over{W-c_i^r}}=
1-(W-c_i^r)^{q-1}.$$  Thus, $$\Delta(W-c_1^r) - \Delta(W-c_2^r) =
(W - c_2^r)^{q-1} - (W - c_1^r)^{q-1}$$
is a polynomial in $Z$ of degree $r(q-2) > 0$ as long as $ q > 2$.

Let $J=(W-c_i^r) A_0$ and $I=M(a_i)^r A_0$.  Then $I^pA_0=J^pA_0$, since after
extension to B, they become equal.  Now $A_0$ is a normal domain, so its group
of invertible fractional ideals is torsion-free.  Hence $I A_0 = J A_0$,
i.e.\ $M(a_i)^rA_0=(W-c_i^r)A_0$.
Therefore, it suffices to show (by \ref{son-of-samuel})
that $\Delta(W-c_1^r) - \Delta(W-c_2^r) \ne \Delta(u)$
for any $u \in A_0^*$.  Since $\Delta(u)$ is a constant
in $Z$, the result follows by the previous paragraph.  \qed
\end{proof}

Now we give an alternate (more geometric) explanation of the preceeding result.
It will show that the kernel of \mapx[[ \Pic(A) || \Pic(A_K) ]] is infinite
and $q$-torsion, but not that the kernel has exponent $q$.

Let $U = \Spec(A)$, and let $V = \Proj(k[X,Y,T]/(X^q-X T^{q-1}-\alpha Y^q))$.
Then $U$ is an open subscheme of $V$.  We have a commutative diagram
\squareSE{\Pic(V)}{\Pic(U)}{\Pic(V_K)}{\Pic(U_K)\nulldot}%
Of course $\Pic(U_K) = 0$ since $A_K$ is a polynomial ring.  Since $U$ is
obtained from $V$ by deleting the single regular point $[X,Y,T] = [1,0,0]$,
\mapx[[ \Pic(V) || \Pic(U) ]] is surjective with cyclic kernel, and in fact one
sees that the kernel is infinite cyclic.  On the other hand, the map
\mapx[[ \Pic(V) || \Pic(V_K) ]] is injective, as is well-known.%
\footnote{More generally, it is even true that if $K/k$ is any field extension,
$V$ is a projective $k$-scheme, and one has two coherent $\O_V$-modules which
become isomorphic over $V_K$, then they were already isomorphic over $V$
-- see \cite{Wi, (2.3)}.}  Hence the kernel of \mapx[[ \Pic(U) || \Pic(U_K) ]]
is exactly $\Pic^0(V)$.

By results of Grothendieck \cite{FGA, (2.1), (3.1)},
there exists a commutative group scheme $P$
of finite type over $k$ such that $\Pic^0(V) = P(k)$.  Since $\dim(V) = 1$,
$H^2(V,\O_V) = 0$, so $P$ is smooth and $\dim(P) = h^1(V,\O_V)$ by
\cite{FGA, \#236, 2.10(ii, iii)}.  Now $h^1(V,\O_V)$ is just the arithmetic
genus of a plane curve of degree $q$, which is $(q-1)(q-2)/2$.  In particular,
since $q \geq 3$, we have $h^1(V,\O_V) > 0$.  Hence $P$ is
positive-dimensional.  Now since $P$ is a geometrically integral scheme of
finite type over a separably closed field $k$, $P(k)$ is Zariski dense in
$P(k^a)$ -- see \cite{Se, discussion on p.\ 107}.  In particular, since $P$ is
positive-dimensional, it follows that $P(k)$ is infinite.  Hence
$\Pic^0(V)$ is infinite.  Hence the kernel of the map
\mapx[[ \Pic(U) || \Pic(U_K) ]] is infinite; it is $q$-torsion by
\pref{kernel-is-torsion}.  This completes the alternate proof of
\pref{deranged-derivations}, except that we have not shown that the kernel
has exponent $q$.  \qed

\vspace{0.03in}

We will see in the next section (at least assuming resolution of singularities)
that for a separated scheme $X$ of finite type over a field $k$, the kernel of
the map $\Pic(X) \to \Pic(X_{k\a})$ is always bounded (i.e., torsion with
finite exponent).  The Picard group itself, however, can have infinite
exponent.  For example, if $X$ is any smooth affine or projective curve
of positive genus over an algebraically closed field, then the torsion
subgroup of $\Pic(X)$ is unbounded, because then $\Pic(X)$ is the quotient of
an abelian variety by a \FG\ subgroup.

\section{Eventual vanishing of the kernel of Pic}\label{eventual-section}

Let $k$ be a field.  The main theorem of this section
\pref{eventually-injective} asserts that if $X$ is a separated $k$-scheme of
finite type, then there exists
a finite field extension $k^+$ of $k$, such that for every algebraic field
extension $L$ of $k^+$, the canonical map $\Pic(X_L) \to \Pic(X_{L\a})$
is injective.

This statement has a sheaf-theoretic formulation, which we consider, in part
because it figures in the proof.  Let $F$ be an
(abelian group)-valued $k$-functor, meaning a functor from
$\cat{$k$-algebras}$ to $\cat{abelian groups}$.
Let $p: B \to C$  be a faithfully flat homomorphism of $k$-algebras.
There are maps $i_1,i_2: C \to C\otimes_BC$, given by
$c \mapsto c \o* 1$ and $c \mapsto 1 \o* c$, respectively.  If the sequence
\diagramno{(*)}{0 & \mapE{} & F(B) & \mapE{F(p)} & F(C) & \mapE{F(i_1)-F(i_2)}
           & F(C\o*_BC)}%
is exact for all $p$, then one says that $F$ is a {\it sheaf\/}
(for the fpqc [faithfully flat quasi-compact] topology).

This is often too much to ask, and so one may look only at certain maps $p$, or
ask only that $F(p)$ be injective.  To relate this to the theorem
\pref{eventually-injective}, let
$F$ be given by $B \mapsto \Pic(X_B)/\Pic(B)$.  Of course if $B$ is a field,
we have $F(B) = \Pic(X_B)$.  What the theorem  says
is that if we enlarge $k$ sufficiently (replacing it by a finite extension),
then $F(L \to L\a)$ is injective, for all algebraic extensions $L$ of $k$.

In this sense, $F$ becomes close to being a sheaf, if we allow for the
enlargement of $k$.  However, in general, one cannot get $(*)$ exact
in an analogous manner.  More precisely, for suitable $k$ and $X$, one cannot
find a finite extension $k^+$ of $k$ such that for every algebraic extension
$L$ of $k^+$, if $p:L \to L\a$ is the canonical map, then $(*)$
is exact.  As an example, let $k$ be a separable closure of $\F_q(t)$,
for some prime $q$, and let $X = \Spec(k[x,y]/(y^2-x^3))$.  Then
$\Pic(X_L) = L$ for every extension field $L$ of $k$.  For any finite extension
$k^+$ of $k$, there is some $a \in k^+ - (k^+)^q$.  If $L = k^+$, one finds
that $a^{1/q}$ lies in the kernel of $F(i_1)-F(i_2)$, but not in the image of
$F(p)$.

\begin{definition}
Let $k$ be a field.  A $k$-scheme $S$ is {\it geometrically stable\/}
if (1) it is of finite type, and (2) every irreducible component of $S\red$ is
geometrically integral and has a rational point.
\end{definition}

\begin{theorem}\label{eventually-injective}
Let $k$ be a field; assume that resolutions of singularities exist for
varieties over $k^a$.  Let $X$ be a separated $k$-scheme of finite type.
Then there exists a finite field extension $k^+$ of
$k$, such that for every algebraic field extension $L$ of $k^+$, the
canonical map $\Pic(X_L) \to \Pic(X_{L\a})$ is injective.
\end{theorem}

\vspace{0.1in}

\par\noindent{\bf Proof of (\ref{eventually-injective})}
In the course of the proof, we will refer to {\it enlarging\/} $k$, by
which we mean that $k$ is to be replaced by a suitably large finite field
extension, contained in $k\a$.  This is done only finitely many times.  Then,
at the end
of the proof, the $k$ we have is really the $k^+$ of which the theorem speaks.
We may treat $L$ as an extension of $k$ which is contained in $k\a$.

By induction on $\dim(X)$, we may assume that the theorem holds
when $X$ is replaced by a scheme of strictly smaller dimension.
(The case of dimension zero is trivial.)

\vspace{0.1in}
\par\noindent{\bf Step 1.  The case of a geometrically normal scheme}
\vspace{0.1in}

If $T$ is a $k$-scheme, we have let $T\a$ denote $T_{k\a}$.
However, in some places in the next paragraph we shall define a $k\a$-scheme
$T\a$, even though $T$ has not yet been defined; we will then
proceed to construct a $k$-scheme $T$ such that $T\a = T_{k\a}$.

Suppose $X$ is geometrically normal.  By \cite{N}, there is a proper
$k\a $-scheme $\oX\a$ which contains $X\a $ as a dense open subscheme.  (Note
that $\oX$ does not yet been defined.)  After normalizing $\oX\a$
\WMAT\ $\oX\a$ is normal.  Let $\pi\a : \oY\a \to \oX\a$ be a
{\it resolution of singularities}, by which we mean that
$\oY\a$ is regular, $\pi\a$ is a proper morphism,  and $\pi\a$ is an
isomorphism over $\oX\a - \Sing(\oX\a)$.  (Note that $\oY$ and $\pi$ have not
yet been defined.)

By looking at the equations defining $\oX\a$, $\oY\a$, and
$\pi\a$, we can (after enlarging $k$ if necessary) find a $k$-scheme
$\oX$
(containing $X$ as a dense open subscheme), a $k$-scheme $\oY$,
and a morphism
$\pi: \oY \to \oX$ such that $\pi \times_k k\a = \pi\a$.
Then $\oX$ is geometrically normal, $\oY$ is
geometrically regular, and (by faithfully flat descent
\cite{EGA4, (2.7.1)(vii)}) $\pi$ is proper.  Let $Y = \pi^{-1}(X)$.

By enlarging $k$ if necessary, we may assume that if $C_1,\dots,C_n$ are the
irreducible codimension one components of $\oY-Y$, and if $p: \oY\a
\to \oY$ is the natural map, then $p^{-1}(C_1), \dots, p^{-1}(C_n)$ are the
irreducible codimension one components of $\oY\a  - Y\a $.

Let $d : \Pic(\oY\a ) \to \Pic(Y\a )$,
$e : \Pic(\oY) \to \Pic(\oY\a )$, $f : \Pic(\oY) \to \Pic(Y)$
and $h : \Pic(Y) \to \Pic(Y\a ) $ be the canonical maps.

For any normal proper $k$-scheme $V$ of finite type, the canonical
map $\Pic(V) \to \Pic(V\a ) $ is injective.  (See
\cite{Mi3, (6.2)}.)  In particular, $e$ is injective.  Since
$$[p^{-1}(C_1)],\ldots,[p^{-1}(C_n)]$$%
generate $\Ker(d)$, it follows
that $[C_1], \ldots, [C_n]$ generate $\Ker(de)$.  Since $f$ is surjective,
it follows that $h$ is injective.

A slight modification of the argument above shows that the canonical
map $h_L : \Pic(Y_L) \to \Pic(Y\a ) $ is injective for every algebraic
extension $L/k$.  Let $\L$ be a line
bundle on $X_L$ that becomes trivial on $X\a $.  Since $h_L$ is
injective,
$\L$ becomes trivial on pullback to $Y_L$.  Let $r : X\a  \to X_L $
be the canonical map.  Then the restriction of
$\L$ to $X_L - r(\Sing(X\a ))$ is trivial.  By (\ref{S2})
$\L$ is trivial.  This completes the proof when $X$ is geometrically
normal.

\vspace{0.1in}
\par\noindent{\bf Step 2.  The case of a geometrically reduced scheme}
\vspace{0.1in}

If $X$ is geometrically reduced we may assume, by enlarging $k$ if
need be, that the normalization
$X\nor$ is geometrically normal.  Let $\pi : X\nor  \to X $ be
the canonical map.  Let $\I = [\O_X : \pi_*\O_{X\nor}]$ be the
conductor of $X\nor $ into $X$.   This
is a coherent sheaf of ideals in
$\O_X$.  Let $X/\I$ denote $\hbox{\bf Spec}(\O_X/\I)$, and let
$X\nor/\I$ denote $(X/\I) \times_X X\nor $.  By further enlarging
$k$, we may assume that  the pullback of $\I$ to $X\a $ is the
conductor of $(X\nor  )\a $
into $X\a $.  Enlarging $k$ still further, we may assume that
$X\nor $, $X/\I$, and $X\nor /\I$ are all geometrically stable.

Let $F_1$ and $F_2$ be the (abelian group)-valued $k$-functors defined by
$$
F_1(B) = \frac{\Gamma((X\nor )_B)^*}{\Gamma(X_B)^*} \text{\ \
and\ \ }
F_2(B) = \frac{\Gamma((X\nor /\I)_B)^*}{ \Gamma((X/\I)_B)^*}.
$$
Let $G = F_2/F_1$, the quotient in $\cat{(abelian group)-valued $k$-functors}$.
We have a commutative diagram
$$
\begin{CD}
{}         @.            0         @.           0        @.        {}
 @.         {}\\
@.                       @VVV               @VVV              @.
       @. \\
0       @>>>       F_1(L)   @>>>   F_2(L)  @>>>  G(L)        @>>>      0 \\
@.                       @VVV                @VVV           @VVV
   @. \\
0      @>>>      F_1(L\a) @>>>  F_2(L\a) @>>>  G(L\a)     @>>>     0 \\
 @.                     @VVV                @VVV             @VVV
   @. \\
0   @>>> F_1(L\a\o*_LL\a)  @>>>     F_2(L\a\o*_LL\a) @>>>
               G(L\a\o*_LL\a)   @>>> 0
\end{CD}
$$
with exact rows.  Now we use \cite{J2, (4.5)}:

\vspace{0.1in}
\par\noindent{\bf Theorem}\ \ {\it
Let $S$ and $T$ be geometrically stable $k$-schemes, and let $f: S \to T$ be a
dominant morphism of $k$-schemes.  Let $Q$ be the (abelian group)-valued
$k$-functor given by $Q(A) = \Gamma(S_A)^* / \Gamma(T_A)^*$.  Let $p: B \to C$
be a faithfully flat homomorphism of reduced $k$-algebras.  Then the sequence
$$0 \to Q(B) \to Q(C) \to Q(C \otimes_B C)$$%
is exact.}
\vspace{0.1in}

This implies that the first two columns are exact.  It follows that the
canonical map $G(L) \to G(L\a)$ is injective.

We need a scheme-theoretic version of  Milnor's Mayer-Vietoris sequence
\cite{Bas1, Chap.\ IX, (5.3)}, which may be found in
\cite{We, (7.8)(i)}, and which implies that there is an exact sequence
$$0 \to F_1(L) \to F_2(L) \to \Pic(X_L) \to
\Pic((X\nor)_L) \times \Pic((X/\I)_L) \eqno(\diamondsuit)$$\label{diamond}%
for each algebraic extension $L/k$.   By the induction
hypothesis announced near the beginning of the proof,
\WMAT\ $\Pic((X/\I)_L) \to \Pic((X/\I)\a )$ is injective.  Also, by
Step 1, $\Pic((X\nor)_L) \to \Pic((X\nor)\a)$ is injective.
Since $G(L) \to G(L\a)$ is injective, it follows that
$\Pic(X_L) \to \Pic(X\a)$ is injective.  Thus the theorem
is true for geometrically reduced schemes.

\vspace{0.1in}
\par\noindent{\bf Step 3. Deal with the nonreduced case}
\vspace{0.1in}

This step may of course be ignored if $X$ is affine.  Otherwise, we
use the following result \cite{J2,(5.2)(a)}:

\vspace{0.1in}
\par\noindent{\bf Theorem}\ \ {\it Let $k$ be a field and $X$ a geometrically
stable $k$-scheme.  Let $i: X_0 \to X$ be a nilimmersion.  For any
$k$-algebra $A$, let $\kappa(A)$ be the kernel of the natural map
$\Pic(X\times_kA) \to \Pic(X_0\times_kA)$.  Let $A \to B$
be a faithfully flat homomorphism of reduced $k$-algebras.  Then the
induced map $\kappa(A) \to \kappa(B)$ is injective.}
\vspace{0.1in}

To complete the proof of (\ref{eventually-injective}), we may assume, by
enlarging   $k$ if need be, that  $X$ is geometrically stable and that $k$
is big enough so that the conclusion is valid for the (geometrically
reduced) scheme $X\red$.  We have the following commutative
diagram, for any algebraic extension $L/k$:
$$
\begin{CD}
\Pic(X_L) @>\alpha>> \Pic((X\red)_L)\\
@V\gamma VV              @V\delta VV\\
\Pic(X\a)  @>\beta>>   \Pic(X\red\a)
\end{CD}
$$
Taking $(A \to B) = (L \to k\a)$ in \cite{J2,(5.2)(a)}, cited above, we see
that $\gamma$ is one-to-one on $\Ker(\alpha)$.  Since
$\delta$ is injective, so is $\gamma$. \qed

\begin{remark}\label{woofwoof}
Let $F$ be a finite field of characteristic different from $5$ and containing
a primitive fifth root of unity $\zeta$.  Let $Y\subset {\Bbb P}^3_F$
be the Fermat quintic given by the equation:
$$x^5+y^5+z^5+t^5=0.$$%
Then the group of fifth roots of unity acts on $Y$ by sending $(x,y,z,t)$
to $(x,\zeta y,\zeta ^2z,\zeta^3t)$.  This action has no fixed points, and
the quotient $X$ is a smooth projective surface which is called the
{\it Godeaux surface}.    For any finite extension $E/F$,
$\CH^2(X_E)_{\op tors}=\Z/5\Z$ \cite{KS1, Proposition 9}.  On the other hand,
over an algebraic closure ${\overline F}$, we have
$\CH^2(X_{\overline F})_{\op tors})=0$ \cite{Mi2}.  Thus
for each $E$ there exists a finite extension $H$ over which the $\Z/5\Z$
dies, but a new one comes to take its place.  Now for any smooth surface over a
field, the natural map $K_0(X)\to \CH^*(X)$
is an isomorphism between the Grothendieck ring and the Chow ring.  Hence
we produce an element of $K_0(X)$ with similar properties.  Thus the theorem
does not hold with $K_0$ in place of $\Pic$.
\end{remark}

\begin{corollary}\label{bounded}
Let $k$ be a field; assume that resolutions of singularities exist for
varieties over $k^a$.  Let $X$ be a $k$-scheme of finite type.
Then the kernel of the map $\Pic(X) \to \Pic(X^a)$ is a bounded torsion group.
\end{corollary}

\begin{proof}
One can (details omitted) use \pref{open-cover} to reduce to the case where
$X$ is separated.  Let $K = k\a$.
Let $k^+$ be as in \pref{eventually-injective}.  Then $\Pic(X) \to
\Pic(X_K)$ and $\Pic(X) \to \Pic(X_{k^+})$ have the same kernel.   By
(\ref{kernel-is-torsion}) the kernel has exponent dividing $[k^+ : k]$. \qed
\end{proof}

\section{Finitely Generated Fields}\label{abs-section}

We say a field $k$ is {\it absolutely finitely generated\/}
if it is finitely generated over its prime subfield.  In this section we will
study the structure of $\Pic(X)$, where $X$ is a scheme of finite type over an
absolutely finitely generated field.  We begin with a result
that is presumably well known, but for which we have found
no reference.

\begin{proposition}\label{normal-over}
Let $X$ be a normal scheme of finite
type over $\Z$ or over an absolutely finitely generated
field $k$. Then $\Pic(X)$ is finitely generated.
\end{proposition}

\begin{proof}
A theorem due to Roquette
\cite{Ro}, \cite{L, Chap. 2, (7.6)} handles the case of
schemes of finite type over $\Z$.  Suppose now that
$X$ is of finite type over the absolutely finitely generated
field $k$.
By \pref{S2} the map $\Pic(X) \to \Pic(X - \Sing(X))$ is
injective.  Therefore \WMAT\ $X$ is regular.  There exists
a finitely generated $\Z$-algebra $A \subset k$ and an
$A$-scheme $X_0$ of finite type such that $X\cong X_0\times_A k$.
We have $X \subset X_0$.  Since $X$ is regular, $\O_{X_0,x}$ is
regular for every $x\in X$.  Since $X_0$ is excellent, its
regular locus is open, so there exists an open subscheme of $X_0$
which is regular and contains $X$.  By replacing $X_0$ by this
subscheme, we may assume that $X_0$ is regular.  The map on
divisor class groups $\Cl(X_0) \to \Cl(X)$ is certainly surjective,
and since both $X_0$ and $X$ are regular, the map
$\Pic(X_0) \to \Pic(X)$ is surjective.
Since $\Pic(X_0)$ is finitely generated (by Roquette's theorem),
so is $\Pic(X)$.  \qed
\end{proof}

The following examples show, in contrast, that the torsion
subgroup of $\Pic(X)$ need not be finite if $X$ is not normal.

\begin{example}
Let $B = \Q[x,x^{-1}]$, and put $A = \Q+(x-1)^2B$.  Then $A$ is a
one-dimensional domain, finitely generated as a $\Q$-algebra, and
$\Pic(A) \cong \Q/\Z$.  In particular, $\Pic(A)$ is an unbounded torsion group.
\end{example}

\begin{proof}
We note that $I:= (x-1)^2B$ is the conductor of $B$ into $A$.  Therefore by
Milnor's Mayer-Vietoris exact sequence \cite{Bas1, Chap. IX, (5.3)}
[or see ($\diamondsuit$, p.\ \pageref{diamond}) for the scheme-theoretic
version], we have $\Pic(A) \cong (B/I)^*/((A/I)^*U)$,
where $U$ is the image of $B^*$ in $(B/I)^*$. But $(A/I)^* =
\Q^*$, so $\Pic(A) \cong (B/I)^*/U$.  Now
$$B^* = \{sx^j:s\in\Q^*,j\in \Z\}\cong \Q^*\oplus \Z,$$%
and $(B/I)^* \cong \Q^*\oplus W$, where $W =\{1+s(x-1):s\in\Q\}
\cong \Q$. By keeping track of these identifications, one
easily gets $\Pic(A) \cong \Q/\Z$.  \qed
\end{proof}

By a slight modification we get an example of finite type over $\Z$:

\begin{example}
Fix a positive integer $m$, put
$B = \Z[x,x^{-1},\frac{1}{m}]$, and let $A = \Z[\frac{1}{m}]+(x-1)^2B$.
Then $A$ is a two-dimensional domain finitely generated as a $\Z$-algebra, and
$\Pic(A) \cong \bigoplus_{p|m}\Z_{p^{\infty}}$.
\end{example}

The pathology in the examples above stems from the fact that $B/I$ is not
reduced.

Before stating our main finiteness theorems
[\pref{seminormal-S2} and \pref{absolutely-p}] we record the following result
from \cite{CGW, (7.4)}:

\begin{theorem}\label{CGW}
Let $k$ be an absolutely \FG\ field and let $\Lambda$ be a finite-dimensional
reduced $k$-algebra.  Let $E_1$ and $E_2$ be intermediate
subalgebras of $\Lambda/k$.
\begin{enumerate}
\item If $k$ has positive characteristic $p$, then $\Lambda^*/E_1^*E_2^*$
is a direct sum of a countably generated free abelian group, a
finite group, and a bounded $p$-group.
\item If $\Lambda/k$ is separable, then $\Lambda^*/E_1^*E_2^*$ is
a direct sum of a countably generated free abelian group and a finite group.
\end{enumerate}
\end{theorem}

\begin{theorem}\label{seminormal-S2}
Let $k$ be a field finitely generated over $\Q$ and let $X$ be a reduced
$k$-scheme of finite type
which is seminormal and $S_2$.  Then $\Pic(X)$ is isomorphic to the direct sum
of a free abelian group and a finite abelian group.
\end{theorem}

\begin{proof}
Let $\I$ be the conductor of $X\nor$ into $X$
(see \S\ref{eventual-section}, Step 2).
Let $X/\I := \hbox{\bf Spec}(\O_X/\I)$ and
$X\nor/\I := (X/\I) \times_X X\nor$ denote the
corresponding closed subschemes.  Since $X$ is seminormal, it follows
\cite{T, (1.3)} that $X/\I$ is reduced.
Let $Q$ be the non-normal locus of $X/\I$, which has codimension $\geq 2$
in $X$.  By \pref{S2}, the canonical map $\Pic(X) \to \Pic(X-Q)$ is injective.
Therefore we may replace $X$ by $X-Q$ and start the proof over, with the
added assumption that $X/\I$ is normal.

Let $D$ be the kernel of the map
$\phi: \Pic(X) \to \Pic(X/\I) \times \Pic(X\nor)$.
By (\ref{normal-over}), the target of this morphism is finitely generated.
Let $\Lambda$ be the integral closure of $k$ in $\Gamma(X\nor/\I)$,
$E_1$  the integral closure of $k$ in
$\Gamma(X/\I)$, and $E_2$ the image in $\Lambda$ of the
integral closure of $k$ in $\Gamma(X\nor)$.  Using \pref{roquette}
and the exact sequence ($\diamondsuit$, p.\ \pageref{diamond}) with $L = k$,
we see that there is an exact sequence
$$\hbox{finitely generated} \to \Lambda^*/E_1^*E_2^* \to D
    \to \hbox{finitely generated}.$$%
By (\ref{CGW})(2), $\Lambda^*/E_1^*E_2^*$ is free $\oplus$ finite.  It follows
that
$D$ and thence $\Pic(X)$ is free $\oplus$ finite.  \qed
\end{proof}

\begin{theorem}\label{absolutely-p}
Let $k$ be an absolutely finitely generated
field of positive characteristic $p$, and let $X$ be a $k$-scheme
of finite type.  Then $\Pic(X)$ has the form
$$\hbox{(countably generated free abelian group)} \oplus
  \hbox{(bounded $p$-group)} \oplus \hbox{(finite group)}.$$
\end{theorem}

\begin{proof}
Let $\cal C$ be the class of abelian groups having the form ascribed to
$\Pic(X)$ in the theorem.  We leave to the reader to verify that $\cal C$
is closed under formation of subgroups and extensions.

Induct on $\dim(X)$; the case where $\dim(X) = 0$ is trivial.  We will reduce
to the case where $X$ is reduced.  For this, it is enough to show $(*)$ that
if $\cal J \IN \O_X$ is a square-zero ideal,
$X_0 = \hbox{\bf Spec}(\O_X/\cal J)$, and $\Pic(X_0) \in \cal C$, then
$\Pic(X) \in \cal C$.  The standard exact sequence of sheaves on $X$
$$\begin{CD}
0 @>>> \J @>{a\kern2pt \mapsto 1+a}>> \O_X^* @>>> (\O_X)/\J @>>> 1
\end{CD}$$
yields on taking cohomology an exact sequence
$$H^1(X,\J) \to H^1(X,\O_X^*) \to H^1(X,(\O_X/\J)^*),$$%
from which $(*)$ follows, since $H^1(X,\J)$ is an $\F_p$-vector space.
Therefore we may assume that $X$ is reduced.

Let $\I, D, \phi$, etc.\ be as in the proof of \pref{seminormal-S2}.
(Here we do not know that $X/\I$ is normal.)  By induction and
\pref{normal-over}, the target of $\phi$ is in $\cal C$, and therefore
$\Im(\phi)$ is also.  Since $\cal C$ is closed under extensions, it suffices
to show that $D \in \cal C$.  Arguing as in the proof of
\pref{seminormal-S2}, with \pref{CGW}(2) replaced by \pref{CGW}(1), we see that
this is the case.  \qed
\end{proof}

\begin{remark}
For $k = \F_p$, it was shown in \cite{J1, (10.11)} that $\Pic(X)$ has the
form
$$\left( \oplus_{n=1}^\infty F \right) \oplus \hbox{(\FG\ abelian group)},$$
where $F$ is a finite $p$-group.
\end{remark}

\section{Complements on $K_0(X)$}

Let $k$ be a field, and let $X$ be a $k$-scheme of finite type.
Let $K_0(X)$ denote the Grothendieck group of vector bundles on $X$.
In this section, which is purely expository, we consider the analog for
$K_0(X)$ of the absolute finiteness results for $\Pic(X)$ proved in sections
\ref{torsion-section} and \ref{abs-section}.  Consider the following table:
\vspace{0.15in}

%\def\scaletopage#1{\par\noindent{\psset{unit=\textwidth}\scaleboxto(1,0){#1}}}
%\def\phraseat(#1,#2)#3{\rput(#1,#2){\parbox{3.3cm}{\centering\small #3}}}
%\def\phraseatl(#1,#2)#3{\rput(#1,#2){\parbox{6cm}{\centering\small #3}}}
%\scaletopage{%
%\psset{xunit=1.04pt, yunit=0.88pt}
%\begin{pspicture}(425,200)
%\psset{linewidth=.3pt}
%\psline(0,0)(425,0)         \psline(0,60)(425,60)     \psline(0,140)(425,140)
%\psline(50,200)(425,200)   \psline(0,0)(0,140)       \psline(50,0)(50,200)
%\psline(150,0)(150,200)     \psline(240,60)(240,200)
%%\psline(335,60)(335,200)
%\psline(425,0)(425,200)
%\psset{linewidth=1.2pt}       \psline(50,0)(425,0)(425,140)(50,140)(50,0)
%\phraseat(25,100){$X$\\ arbitrary}  \phraseat(25,30){$X$\\ regular}
%\phraseat(100,170){$k$ algebraically\\ closed}
%\phraseat(195,170){$k$ absolutely\\ \FG\\ of char.\ $p > 0$}
%\phraseat(287.5,170){$k$ finite}
%\phraseat(380,170){$k$ absolutely\\ \FG\\ of characteristic $0$}
%\phraseatl(287.5,30){\circno5 finitely generated}
%\phraseat(100,100){\circno1 finite $n$-torsion for all $n$\\ invertible in
%%$k$}
%\phraseat(195,100){\circno3 finite torsion mod $p^n$-torsion,\\ for some $n$}
%\phraseat(287.5,100){\circno4 \FG\ mod $p^n$-torsion,\\ for some $n$}
%\phraseat(100,30){\circno2 finite $n$-torsion for all $n$}
%\end{pspicture}}

\vspace{0.05in}

% Delete the following table if you restore the pstricks code.
\begin{verbatim}
              k alg.            k abs. f.g.        k finite       k abs f.g.
              closed            of char. p>0                      of char. 0

X arbitrary   [1] finite        [3] finite         [4] f.g.
              n-torsion for     torsion mod        mod p^n-torsion
              all n             p^n-torsion,       for some n
              invertible in k   for some n

X regular     [2] finite
              n-torsion         [5]  f i n i t e l y      g e n e r a t e d
              for all n
\end{verbatim}

If we view this as a collection of statements about $\Pic(X)$, then all
five statements are true\footnote{For statement \circno2, we have used
resolution of singularities.}, as we have seen in the preceeding sections.
{\bf From now on, regard the table as a table of five conjectures about
$K_0(X)$.}  The numbering of these conjectures is not related to the numbering
of results in the introduction.

It would be surprising if all five of these conjectures held.  There is no
field over which any of them are known to hold, even if one restricts
attention to smooth projective or smooth affine schemes $X$.  Conjecture
\circno5\ would follow from a conjecture of Bass \cite{Bas2, \S9.1}
to the effect that $K_i(X)$ is \FG\ for
all $i \geq 0$ and all regular schemes $X$ which are of finite type over $\Z$.

Over an arbitrary field $k$, it is not clear what sort of finiteness statement
might hold for $K_0(X)$: there are examples of infinite $n$-torsion in the
Chow groups of smooth projective varieties over a field of characteristic
zero \cite{KM}.

If $X$ is smooth of dimension $n$, then the operation of taking Chern classes
defines a homomorphism of graded rings \mapx[[ \Gr[K_0(X)] || \CH^*(X) ]],
which becomes an isomorphism after tensoring by $\Z[1 / (n-1)!]$.
(See \cite{F, (15.3.6)}.)  It follows that for any almost
any question about $K_0(X)$, there is a parallel question about the groups
$\CH^q(X)$, $q = 1, \ldots, n$.  Note that when $q = 1$, $\CH^1(X) = \Pic(X)$.

For the remainder of this section, suppose that $X$ is smooth, projective,
and of dimension $n$; we give a partial discussion of results and conjectures
about $\CH^q(X)$, for $q \geq 2$.

First suppose that $k$ is algebraically closed.  Some things are known when
$q \in \setof{2,n}$:  (i) The group ${}_m \CH^2(X)$ is finite if $m$ is
invertible in $k$ \cite{Ra, (3.1)}; (ii) the group
${}_m \CH^n(X)$ is finite for every $m$.  This follows from Roitman's
theorem (See e.g.{\ } \cite{Ra, (3.2)}.)  Hence
\circno2\ holds for smooth projective surfaces.

Now suppose that $k$ is a number field.  Bloch has conjectured that
$\CH^q(X)$ is \FG\ for
all $q$.  All results in this direction assume at least that $H^2(X,\O_X) = 0$.
With this hypothesis, it has been shown that the torsion subgroup of $\CH^2(X)$
is finite \cite{CR}.  Moreover, if $X$ is a surface which is not
of general type, it is known \cite{Sal}, \cite{CR} that $\CH^2(X)$ is \FG.
Hence conjecture \circno5\ holds for a smooth projective surface over a
number field which is not of general type.

Finally, suppose that $k$ is a finite field.  Again, it is conjectured that
$\CH^q(X)$ is \FG\ for all $q$.  What is known is that $\CH^n(X)$ is
\FG\ (in fact $\CH_0(X)$ is \FG\ for any scheme $X$ of finite type over $\Z$
\cite{KS2}), and that the torsion subgroup of $\CH^2(X)$ is finite \cite{CSS},
\ cf.{\ }\cite{CR, (3.7)}.  The first assertion implies that conjecture
\circno5\ holds for any smooth projective surface over a finite field.

\end{document}